\documentclass[sigconf]{acmart}

\setcopyright{none}
\acmYear{2021} 
\copyrightyear{The DOI of the original article is as follows} 
\acmConference[This is the author version of the paper accepted for GECCO '21]{2021 Genetic and Evolutionary Computation Conference}{July 10--14, 2021}{Lille, France}
\acmPrice{}
\acmDOI{10.1145/3449639.3459266}
\acmISBN{}

\usepackage{balance}

\usepackage{graphicx}
\graphicspath{{../}} 

\usepackage{subcaption}
\usepackage{booktabs} 
\usepackage{longtable}

\usepackage{latexsym}
\usepackage{amsmath} 
\usepackage{mathtools}       
\usepackage{stmaryrd}        
\usepackage{bm}              

\providecommand{\argmin}{\operatornamewithlimits{argmin}} 
\providecommand{\R}{\mathbb{R}} 
\providecommand{\E}{\mathbb{E}} 
\providecommand{\T}{\mathrm{T}} 
\renewcommand{\geq}{\geqslant} 
\renewcommand{\leq}{\leqslant} 
\DeclarePairedDelimiterX{\inner}[2]{\langle}{\rangle}{#1, #2}
\DeclarePairedDelimiter{\norm}{\lVert}{\rVert}
\DeclarePairedDelimiter{\abs}{\lvert}{\rvert}

\usepackage{xcolor}
\usepackage{ifthen}

\usepackage{cleveref}

\usepackage{algorithm}
\usepackage{algpseudocode}

\begin{document}

\title{Saddle Point Optimization with Approximate Minimization Oracle}

\author{Youhei Akimoto}
\affiliation{%
  \institution{University of Tsukuba \& RIKEN AIP}
  \streetaddress{1-1-1 Tennodai}
  \city{Tsukuba} 
  \state{Ibaraki, Japan} 
  \postcode{305-8573}
}
\email{akimoto@cs.tsukuba.ac.jp}


\begin{abstract}
A major approach to saddle point optimization $\min_x\max_y f(x, y)$ is a gradient based approach as is popularized by generative adversarial networks (GANs). In contrast, we analyze an alternative approach relying only on an oracle that solves a minimization problem approximately. Our approach locates approximate solutions $x'$ and $y'$ to $\min_{x'}f(x', y)$ and $\max_{y'}f(x, y')$ at a given point $(x, y)$ and updates $(x, y)$ toward these approximate solutions $(x', y')$ with a learning rate $\eta$. On locally strong convex--concave smooth functions, we derive conditions on $\eta$ to exhibit linear convergence to a local saddle point, which reveals a possible shortcoming of recently developed robust adversarial reinforcement learning algorithms. We develop a heuristic approach to adapt $\eta$ derivative-free and implement zero-order and first-order minimization algorithms. Numerical experiments are conducted to show the tightness of the theoretical results as well as the usefulness of the $\eta$ adaptation mechanism.
\end{abstract}

%
%

\ccsdesc[500]{Mathematics of computing~Continuous optimization}
\ccsdesc[500]{Theory of computation~Convergence and learning in games}
\ccsdesc[500]{Theory of computation~Theory of randomized search heuristics}

\keywords{Minimax Optimization, Saddle Point Optimization, Robust Optimization, Convergence}

\maketitle

\providecommand{\X}{\R^m}
\providecommand{\Y}{\R^n}
\providecommand{\Z}{\R^{\ell}}

\sloppy

\section{Introduction}

We consider the following min--max optimization problem
\begin{equation}
  \min_{x \in \X} \max_{y \in \Y} f(x, y) \enspace,\label{eq:spo}
\end{equation}
where $f: \X \times \Y \to \R$ is the objective function.
This arises in many fields, including constrained optimization \cite{Cherukuri2017sicon}, robust optimization \cite{Conn2012bilevel,mmde2018}, robust reinforcement learning (RL) \cite{Pinto2017icml,Shioya2018iclr},
and generative adversarial networks (GANs) \cite{goodfellow2014generative,Salimans2016nips}.

Unless $f$ is convex--concave, i.e., convex in $x$ and concave in $y$, it is often computationally intractable to locate the global optimal solution.
A more realistic target is to locate a local min--max saddle point $(x^*, y^*)$, which is a local minimum of $f$ in $x$ and a local maximum of $f$ in $y$.
In this study, we focus on locating a local min--max saddle point of \eqref{eq:spo}.
Unless otherwise specified, by ``a saddle point,'' we mean a min--max saddle point in this paper. 

First-order approaches, as popularized by the success of GANs, are often employed for this purpose. A simultaneous gradient method
\begin{equation}
  (x_{t+1}, y_{t+1}) = (x_t, y_t) + \eta (- \nabla_x f(x_t, y_t), \nabla_y f(x_t, y_t)), 
  \label{eq:grad}
\end{equation}
has been analyzed for its local and global convergence properties.
On $C^2$ functions, the continuous time dynamics obtained for the limit $\eta \to 0$ are analyzed and their local asymptotic stability around strict local saddle points is shown \cite{Nagarajan2017nips,Cherukuri2017sicon}.
The analysis is extended to \eqref{eq:grad} with a finite $\eta > 0$, and the condition on $\eta$ for a strict local saddle point to be asymptotically stable is derived in \cite{Mescheder2017nips}. 
\citet{pmlr-v89-adolphs19a} has shown that there are locally asymptotically stable points that are not saddle points and proposed a modification to avoid these undesired stable points.
\citet{liang2019interaction} have shown that \eqref{eq:grad} globally converges toward the strict saddle point if the objective function is globally and strongly convex--concave and that some improved gradient-based approaches \cite{daskalakis2017training,yadav2017stabilizing,Mescheder2017nips} can converge toward a nonstrict saddle point on a bilinear function.
For constrained cases, the problems are often treated as variational inequalities, and first-order approaches are employed \cite{gidel:hal-01403348,NIPS2019_9631}.

Zero-order approaches for \eqref{eq:spo} are categorized as coevolutionary approaches \cite{Al-Dujaili2019lego,mmde2018,Jensen2004book,Branke2008ppsn}, Bayesian optimization approaches \cite{Picheny2019bo,Bogunovic2018nips}, trust-region approaches \cite{Conn2012bilevel}, and gradient approximation approaches \cite{Liu2020icml}.
They are often designed heuristically, and little attention has been paid to convergence guarantees and convergence rate analysis of these approaches.
It is known that coevolutionary approaches suffer from nonconvergent behavior \cite{Al-Dujaili2019lego,mmde2018}.
Recently, \cite{Bogunovic2018nips} showed regret bounds for a Bayesian optimization approach and 
\cite{Liu2020icml} showed an error bound for a gradient approximation approach, where the error is measured by the square norm of the gradient. 
Both analyses show sublinear rates under possibly stochastic (i.e., noisy) versions of \eqref{eq:spo}.

In contrast to previous works, we analyze an alternative approach that relies solely on approximate minimization oracles, which can be zero-order, first-order, or higher-order approaches.
Given $(x, y)$, this approach first locates approximate solutions $x'$ and $y'$ to the minimization problems $\argmin_{x' \in \X} f(x', y)$ and $\argmin_{y' \in \Y} - f(x, y')$ by using approximate minimization oracles.
Then, it updates $(x, y)$ toward $(x', y')$ with a learning rate $\eta > 0$, i.e., $(x,y) \leftarrow (1-\eta) (x, y) + \eta (x', y')$.
One can choose domain-specific methods to minimize $f(\cdot, y)$ and $-f(x, \cdot)$.
This work is motivated by a recent advance in robust RL \cite{Pinto2017icml,Shioya2018iclr}, where the objective is to locate a robust policy $x$ under adversarial disturbance $y$, and these parameters are updated in the abovementioned manner; the approximate minimization oracles are some standard RL approaches.
The applications of our analysis are of course not limited to robust RL.
In particular, it is suitable for problems involving numerical simulations and requiring zero-order approaches not necessarily connected to gradient descent.

The contributions of this paper are twofold.
We analyze the local and global convergence properties of an oracle-based saddle point optimization on locally and globally strong convex--concave smooth functions.
We derive a sufficient condition on the learning rate $\eta$ to guarantee linear (i.e., geometric) convergence and derive an upper bound on the convergence rate.
  In contrast to the abovementioned analysis for zero-order approaches showing a sublinear decrease \cite{Bogunovic2018nips,Liu2020icml}, our analysis is for linear convergence; hence, the result is more related to the one obtained for a simultaneous gradient update \cite{Mescheder2017nips,liang2019interaction}.
  We show the condition on $\eta$ to be not only sufficient but also necessary for convergence on a convex--concave quadratic function.
  This reveals a possible shortcoming of approaches with $\eta = 1$, which do not guarantee convergence but are employed in e.g.\ robust adversarial RL \cite{Pinto2017icml,Shioya2018iclr} or coevolutionary approaches.
  The tightness and possible room for improvement are demonstrated in numerical experiments.

  We propose a heuristic method to adapt the learning rate $\eta$. Our approach is derivative free and black box, that is, with no gradient information and no characteristic constant such as a smoothness parameter.
  Because $f$ is often black-box and no sufficient information is available to select the right $\eta$ in advance, especially when zero-order approaches are desired, $\eta$ adaptation is essential for practical approaches.
  We instantiate the whole algorithm using a zero-order randomized hill-climbing approach, namely, the (1+1)-ES with 1/5 success rule \cite{Devroye:72,Schumer:Steiglitz:68,rechenberg:1973} as well as sequential least squares programming (SLSQP) \cite{kraft1988software}.
  We demonstrate that our learning rate adaptation method locates a nearly optimal saddle point by iterating up to three times more $f$-calls compared with the best fixed $\eta$.

\section{Formulation}

\paragraph{Notation}

Suppose that 
$f: \X\times\Y \to \R$ is a twice continuously differentiable function, i.e., $f \in \mathcal{C}^2(\X\times\Y, \R)$. 
Let $H_{x,x}(x,y)$, $H_{x,y}(x,y)$, $H_{y,x}(x,y)$, and $H_{y,y}(x,y)$ be the blocks of the Hessian matrix of $f$, whose $(i, j)$-th component is $\partial^2 f / \partial x_i \partial x_j$, $\partial^2 f / \partial x_i \partial y_j$, $\partial^2 f / \partial y_i \partial x_j$, and $\partial^2 f / \partial y_i \partial y_j$, respectively, evaluated at a given point $(x, y)$.
Let $J_g(z)$ be the Jacobian of a differentiable $g = (g_1, \dots, g_k): \Z \to \R^k$, where the ($i$, $j$)-th element is $\partial g_i / \partial z_j$ evaluated at $z = (z_1,\dots,z_\ell) \in \Z$.
If $k = 1$, we write $J_g(z) = \nabla g(z)^\T$.
For a positive definite symmetric matrix $A$, let $\sqrt{A}$ denote the matrix square root.

\paragraph{Saddle Point}

A point $(x^*, y^*) \in \X \times \Y$ is a (strict) local saddle point of a function $f: \X \times \Y \to \R$ if there exists a neighborhood $\mathcal{E}_x \times \mathcal{E}_y \subseteq \X \times \Y$ including $(x^*, y^*)$ such that for any $(x, y) \in \mathcal{E}_x\times \mathcal{E}_y$, the condition $f(x, y^*) \geq (>) f(x^*, y^*) \geq (>) f(x^*, y)$ holds.
If $\mathcal{E}_x = \X$ and $\mathcal{E}_y = \Y$, $(x^*, y^*)$ is called the (strict) global saddle point.
For $f \in \mathcal{C}^2(\X\times\Y, \R)$, a point $(x^*, y^*)$ is a strict saddle point if it is a critical point ($\nabla_x f(x^*, y^*) = 0$ and $\nabla_y f(x^*, y^*) = 0$) and $H_{x,x}(x^*, y^*) \succ 0$ and $H_{y,y}(x^*, y^*) \prec 0$ both hold.

\paragraph{Convex--concave Function}

A function $f:\X\times\Y\to\R$ is a (strict) convex--concave function if $f$ is (strictly) convex in $x$ and (strictly) concave in $y$. It is locally (and strictly) convex--concave in open sets $\mathcal{E}_x \subseteq \X$ and $\mathcal{E}_y \subseteq \Y$ if the restriction of $f$ to $\mathcal{E}_x \times \mathcal{E}_y$ is a (strict) convex--concave function. Moreover, $f$ is called a (locally) strong convex--concave function if (the restriction to $\mathcal{E}_x \times \mathcal{E}_y$ of) $f$ is strongly convex in $x$ and strongly concave in $y$.
For $f \in \mathcal{C}^2(\X\times\Y, \R)$, $f$ is (locally) strong convex--concave if and only if there exists $\mu > 0$ such that $H_{x,x}(x, y) \succcurlyeq \mu$ and $H_{y,y}(x, y) \preccurlyeq - \mu$ for all $(x, y) \in \X\times\Y$ ($(x, y) \in \mathcal{E}_x \times \mathcal{E}_y$).
  
\paragraph{Suboptimality Error}

We define a quantity to measure the progress of the optimization towards the saddle point.
For a strict convex--concave function, there exists a unique saddle point $(x^*, y^*)$.
The suboptimality error \cite{gidel:hal-01403348} is then defined as
\begin{align}
  G(x, y) &:= \max_{y' \in \Y} f(x, y') - \min_{x' \in \X} f(x', y)
    \enspace. \label{eq:subopt-global}
\end{align}
It is easy to see that the suboptimality error is non-negative and is zero only at the saddle point $(x^*, y^*)$.
If $f$ is locally strictly convex--concave in a neighborhood $\mathcal{E}_x \times \mathcal{E}_y$ of a local saddle point $(x^*, y^*)$, the suboptimality error can be extended to a nonconvex--concave scenario as
\begin{align}
  G(x, y) &:= \max_{y' \in \mathcal{E}_y} f(x, y') - \min_{x' \in \mathcal{E}_x} f(x', y) 
    \enspace. \label{eq:subopt-local}
\end{align}

In this study, we use a quadratic approximation of \eqref{eq:subopt-global} and \eqref{eq:subopt-local} to measure the progress toward a saddle point $(x^*, y^*)$.
If $f \in \mathcal{C}^2(\X\times\Y, \R)$ and (locally) strong convex--concave around a (local) saddle point $(x^*, y^*)$,
the (local) suboptimality error can be approximated around $(x^*, y^*)$ by
\begin{equation}
    \tilde{G}(x, y) = \norm{x - x^*}_{G_{x,x}^*}^2/2 + \norm{y - y^*}_{G_{y,y}^*}^2/2 \enspace,
    \label{eq:gap-quad}
  \end{equation}
where (dropping $(x, y)$ from $H_{x,x}$, $H_{y,y}$, $H_{x,y}$, and $H_{y,x}$)
\begin{equation}
  \begin{split}
    G_{x,x}(x,y) =& H_{x,x} + H_{x,y}(-H_{y,y})^{-1}H_{y,x}\\
    G_{y,y}(x,y) =& (-H_{y,y}) + H_{y,x}H_{x,x}^{-1}H_{x,y}
  \end{split}\label{eq:gmat}
\end{equation}
and
\begin{equation}
  \begin{split}
    \norm{x - x^*}_{G_{x,x}^*} =& [(x - x^*)^\T G_{x,x}(x^*,y^*) (x - x^*)]^{1/2} ,\\
    \norm{y - y^*}_{G_{y,y}^*} =& [(y - y^*)^\T G_{y,y}(x^*,y^*) (y - y^*)]^{1/2} .
  \end{split}\label{eq:gap-norm}
\end{equation}
It is non-negative for all $(x, y) \in \X\times\Y$ and zero only at $(x^*, y^*)$. 
If $f$ is a convex--concave quadratic function,
we have $G(x, y) = \tilde{G}(x,y)$ for all $(x,y) \in \X\times\Y$.
In the following we use $\tilde{G}(x, y)$ to measure progress.

\section{Algorithm}\label{sec:oracle}

\paragraph{Approximate Minimization Oracle}

Consider a minimization of $g: \Z \to \R$.
Given a positive definite symmetric matrix $A \in \R^{\ell \times \ell}$, consider a local minimization of $g$ in the neighborhood $U_A(z, r) := \{z' \in \Z : \norm{z' - z}_{A} \leq r \}$ of an initial search point $z$ with radius $r > 0$ under $\norm{z}_{A} = [ z^\T A z ]^{1/2}$,
\begin{align}
  z^* = \argmin_{z' \in U_{A}(z,r)} g(z') \enspace.\label{eq:local-opt}
\end{align}
A local minimization oracle $M_{\epsilon,r}^{A}(g, z)$ takes the function $g$ to be minimized, and the initial solution $z$ as inputs and outputs an approximate solution $z'$ to the above minimization problem \eqref{eq:local-opt} satisfying the following condition: 
\begin{align}
  \norm{ M_{\epsilon,r}^{A}(g, z) - z^*}_{A}^2 \leq \epsilon \cdot \norm{z - z^*}_A^2 \enspace.\label{eq:oracle}
\end{align}
In other words, it locates a point that decreases the squared distance from the exact local minimum under $A$ by the factor $\epsilon \in (0, 1)$ compared to the initial point $z$.

Important examples are algorithms that exhibit linear convergence, where the runtime (number of $f$-calls and/or $\nabla f$-calls) to shrink the distance to the optimum by the factor $\epsilon$ is estimated by $\Theta(\log(1/\epsilon))$.
Therefore, by running $\Theta(\log(1/\epsilon))$ iterations of such algorithms, one can implement the above oracle. However, it is not limited to linearly convergent algorithms. The oracle requirement \eqref{eq:oracle} can be satisfied with algorithms that converge sublinearly. In such cases, the number of iterations required to satisfy \eqref{eq:oracle} becomes greater as the candidate solution approaches the optimum. Therefore, the stopping condition ($\tau$ in our proposed algorithm, see \Cref{sec:lra}) for the search algorithm inside the oracle needs to be tuned more carefully.

\paragraph{Oracle-based Saddle Point Optimization}
We consider an approach to the saddle point optimization \eqref{eq:spo} based solely on approximate local minimization oracles satisfying \eqref{eq:oracle}.
We first find approximate local solutions to $\argmin_{x' \in \X} f(x', y_t)$ and $\argmin_{y' \in \Y} - f(x_t, y')$ by approximate minimization oracles $M_{\epsilon_x,r_x}^{A_x}$ and $M_{\epsilon_y,r_y}^{A_y}$, respectively.\footnote{We remark that if the right-hand side of condition~\eqref{eq:oracle} is replaced with a constant, the whole minimax algorithm cannot be guaranteed to converge to a local saddle point. Because our objective is to derive the linear convergence rate, such a situation is beyond the scope of this paper.}
Then, $x_t$ and $y_t$ are updated with the learning rate $\eta > 0$ as 
\begin{equation}
  \begin{split}
  x_{t+1} &= (1 - \eta) x_{t} + \eta \cdot M_{\epsilon_x,r_x}^{A_x}( f(\cdot, y_{t}), x_t) \enspace, \\
  y_{t+1} &= (1 - \eta) y_{t} + \eta \cdot M_{\epsilon_y,r_y}^{A_y}( -f(x_{t}, \cdot), y_t) \enspace.
  \end{split}\label{eq:algo}
\end{equation}




\section{Convergence Analysis}

We investigate the global and local convergence properties of the proposed approach \eqref{eq:algo} on $f \in \mathcal{C}^2(\X\times\Y, \R)$.
We are especially interested in knowing how small the learning rate $\eta$ needs to be and how quickly it converges. 

In the following, let $(x^*, y^*)$ be a strict saddle point of $f$ and $f$ be locally strong convex--concave around $(x^*, y^*)$.
For notational simplicity, let $H_{x,x}^{*} = H_{x,x}(x^*, y^*)$, $H_{y,y}^* = H_{y,y}(x^*, y^*)$, $H_{x,y}^* = H_{x,y}(x^*, y^*)$, $H_{y,x}^* = H_{y,x}(x^*, y^*)$, $G_{x,x}^* = G_{x,x}(x^*, y^*)$, and $G_{y,y}^* = G_{y,y}(x^*, y^*)$. 

\providecommand{\opty}{\hat{y}}
\providecommand{\optx}{\hat{x}}
Our analysis is based on
the implicit function theorem (e.g., Theorem~5 of \cite{deoliveira2013}), which shows the existence and uniqueness of the solutions $\hat{x}(y)$ and $\hat{y}(x)$ to $\nabla_x f(\hat{x}(y), y) = 0$ and $\nabla_y f(x, \hat{y}(x)) = 0$.
The proofs of the following results are provided in \Cref{apdx:proof}. 

\subsection{Global Linear Convergence}

Suppose that $f \in \mathcal{C}^2(\X\times\Y, \R)$ is a strong convex--concave function.
Then, there exists a unique global saddle point $(x^*, y^*)$.
We assume that for
\begin{equation}
  \begin{split}
    \Delta_{x,y}(x,y) &= (H_{x,x}(x, y))^{-1}H_{x,y}(x, y) - (H_{x,x}^*)^{-1}H_{x,y}^* \\
    \Delta_{x,y}(x,y) &= (H_{y,y}(x, y))^{-1}H_{y,x}(x, y) - (H_{y,y}^*)^{-1}H_{y,x}^* \enspace,
  \end{split}
\end{equation}
there exists $\delta \in [0, 1)$ such that 
\begin{equation}
  \begin{aligned}
    \sigma( \sqrt{\smash[b]{G_{x,x}^*}} \Delta_{x,y}(x, y) \sqrt{\smash[b]{G_{y,y}^*}}^{-1} ) \leq \delta\\
    \sigma( \sqrt{\smash[b]{G_{y,y}^*}} \Delta_{y,x}(x, y) \sqrt{\smash[b]{G_{x,x}^*}}^{-1} ) \leq \delta
\end{aligned}
\label{eq:delta}
\end{equation}
 hold for all $(x, y) \in \X\times\Y$, where $\sigma(\cdot)$ denotes the greatest singular value of the argument matrix.
We let the class of such functions be denoted by $\mathcal{F}$. 
An instance in $\mathcal{F}$ is a convex--concave quadratic function, where $\delta = 0$.

The following theorem states the conditions on the approximation error $\epsilon_x$ and $\epsilon_y$ of the approximate minimization oracle \eqref{eq:oracle} and the learning rate $\eta$ to guarantee the global linear convergence and shows the upper bound of the convergence rate.

\begin{theorem}
  \label{thm:glc}
  Suppose $f \in \mathcal{F}$ and $(x^*, y^*)$ is the global saddle point of $f$.
  Let $\tilde{G}(x, y)$ be defined as \eqref{eq:gap-quad}. 
  Let $\bar{\sigma}$ be the greatest singular value of $\sqrt{\smash[b]{G_{x,x}^*}} (H_{x,x}^*)^{-1}H_{x,y}^* \sqrt{\smash[b]{G_{y,y}^*}}^{-1}$, which is equal to $\sqrt{\smash[b]{G_{x,x}^*}}^{-1} H_{x,y}^* (-H_{y,y}^*)^{-1} \sqrt{\smash[b]{G_{y,y}^*}}$.
  Suppose that there are approximate minimization oracles $M_{\epsilon_x,r_x}^{A_x}$ and $M_{\epsilon_y,r_y}^{A_y}$ satisfying oracle condition \eqref{eq:oracle} with $A_x = G_{x,x}^*$ and $A_y = G_{y,y}^*$. 
  Consider algorithm \eqref{eq:algo} with the neighborhood parameter $r_x = r_y = \infty$.
  Let $\bar\epsilon = \max[\epsilon_x, \epsilon_y]$. 
  If 
  \begin{align}\label{eq:glc-condition}
    \bar\epsilon < 1 - \delta &&\mathrm{and}&&
    \eta < \bar\eta := \frac{2 (1 - (\bar\epsilon + \delta))}{ 1+\bar\sigma^2 - (\bar\epsilon + \delta)^2 } \enspace,
  \end{align}
  then for any $(x_0, y_0) \in \X\times\Y$, 
  \begin{equation}\label{eq:glc-conclusion}
    \tilde{G}(x_t, y_t) \leq \gamma^{2t} \tilde{G}(x_0, y_0)
  \end{equation}
  holds for $\gamma = ((1-\eta)^2 + \eta^2 \bar\sigma^2)^{1/2}
  + \eta (\bar\epsilon + \delta) < 1$.
\end{theorem}

\subsection{Local Linear Convergence}

If $f \in \mathcal{C}^2(\X\times\Y, \R)$, and it is locally strong convex--concave in a neighborhood around a local saddle point $(x^*, y^*)$,
we can derive a local convergence condition, as shown in \Cref{thm:llc}.
This provides a tighter estimate of the asymptotic rate $\bar\gamma$ of convergence than that in \Cref{thm:glc}.

\begin{theorem}
  \label{thm:llc}
  Let $f \in \mathcal{C}^2(\X\times\Y, \R)$ and locally strong convex--concave in a neighborhood around a local saddle point $(x^*, y^*)$.
  Let $\tilde{G}(x, y)$ be defined as \eqref{eq:gap-quad}.
  Suppose that there are approximate minimization oracles $M_{\epsilon_x,r_x}^{A_x}$ and $M_{\epsilon_y,r_y}^{A_y}$ satisfying oracle condition \eqref{eq:oracle} with $A_x = G_{x,x}^*$ and $A_y = G_{y,y}^*$.   
  Consider algorithm \eqref{eq:algo} with neighborhood parameters $r_x$ and $r_y$.
  Let $\bar{\sigma}$ and $\bar\epsilon$ be defined in \Cref{thm:glc}. 
  If  
  \begin{align}\label{eq:llc-condition}
    \bar\epsilon < 1
    && \text{and}
    && \eta < \bar\eta := \frac{2 (1 - \bar\epsilon)}{ 1+\bar\sigma^2 - \bar\epsilon^2 } \enspace,
  \end{align}
  then for any $\gamma \in (\bar\gamma, 1)$ with $\bar\gamma = ((1-\eta)^2 + \eta^2 \bar\sigma^2)^{1/2}
  + \eta \bar\epsilon$, 
  there exist constants $\bar{r}_x > 0$ and $\bar{r}_y > 0$ such that
  for any $r_x \in (0, \bar{r}_x)$ and $r_y \in (0, \bar{r}_y)$, there exists a neighborhood $U \subseteq \X\times\Y$ of $(x^*, y^*)$
  satisfying 
  \begin{equation}\label{eq:llc-conclusion}
    \tilde{G}(x_t, y_t) \leq \gamma^{2t} \tilde{G}(x_0, y_0),  \quad \forall (x_0, y_0) \in U. 
  \end{equation}
\end{theorem}

\subsection{Discussion}\label{sec:discussion}

\Cref{thm:glc,thm:llc} show that the approximate suboptimality error $\tilde{G}(x_t, x_t)$ converges linearly in terms of the number of queries to the approximate minimization oracle as long as the conditions are satisfied.
If the approximate minimization oracle requires $\mathcal{O}(1)$ $f$-calls (or $\nabla f$-calls) on average, this implies that $\tilde{G}(x_t, y_t)$ converges linearly in terms of $\# f$-calls (or $\nabla f$-calls) as well.
As mentioned in \Cref{sec:oracle}, a linearly converging minimization approach requires $\mathcal{O}(\log(1/\epsilon))$ $f$-calls, and hence the resulting approach \eqref{eq:algo} converges linearly in terms of queries.
This result is distinguished from the results of \citet{Bogunovic2018nips,Liu2020icml}, where a sublinear decrease has been analyzed in a stochastic (noisy) setting, rather than a linear decrease in a deterministic one. 

\Cref{thm:llc} provides the upper bound $\bar\gamma$ of the convergence rate $\lim_{t\to\infty} \sqrt{\smash[b]{\tilde{G}(x_{t+1}, y_{t+1}) / \tilde{G}(x_t, y_t)}}$.
It is minimized at $\eta = \eta^*$, and its optimal value is $\bar\gamma^*$, where
\begin{equation}
  \begin{split}
  \eta^* &:= \frac{ 1 - \sqrt{ \bar\sigma^2 \bar\epsilon^2 / (1 - \bar\epsilon^2 + \bar\sigma^2)} }{ 1+\bar\sigma^2 } \enspace,\\
  \bar\gamma^* &:= \frac{\bar\sigma(1+\bar\sigma^2) + (\bar\epsilon + \bar\epsilon^2\bar\sigma)\sqrt{1 - \bar\epsilon^2 + \bar\sigma^2}}{(1+\bar\sigma^2) \sqrt{1 - \bar\epsilon^2 + \bar\sigma^2}} \enspace.
  \end{split}\label{eq:optimal}
\end{equation}
We have $\eta^* = 1/(1 + \bar\sigma^2)$ and $\bar\gamma^* = \bar\sigma / \sqrt{1 + \bar\sigma^2}$ for $\bar\epsilon = 0$ (an exact minimization oracle) and $\eta^* \downarrow 0$ and $\bar\gamma^* \uparrow 1$ as $\bar\epsilon \uparrow 1$.
The number of iterations to halve $\tilde{G}$ is $T \leq \log(2) / \log(1/\bar\gamma^2)$.

Conditions \eqref{eq:glc-condition} and \eqref{eq:llc-condition} are shown to be sufficient in \Cref{thm:glc,thm:llc}. 
The question naturally arises as to whether it is necessary as well.
For example, consider $f(x,y) = x^2/2 + \bar\sigma xy - y^2/2$.
Then, $G_{x,x}^* = G_{y,y}^* = 1 + \bar\sigma^2$, where $\bar\sigma$ is equivalent to $\bar\sigma$ in \Cref{thm:glc}.
We have $\hat{y}(x) = \bar\sigma x$ and $\hat{x}(y) = - \bar\sigma y$.
For the exact minimization oracle ($\bar\epsilon = 0$), we have
$x_{t+2} = [(1-\eta)^2 + \eta^2\bar\sigma^2] x_{t}$ and $y_{t+2} = [(1-\eta)^2 + \eta^2\bar\sigma^2] y_{t}$,
and it is easy to see that $(x_{t}, y_{t})$ converges if and only if $(1-\eta)^2 + \eta^2\bar\sigma^2 < 1$, i.e., $\eta < 2 / (1 + \bar\sigma^2)$, the right-hand side of which is $\bar\eta$ in \eqref{eq:glc-condition} and \eqref{eq:llc-condition}. 
Therefore, conditions \eqref{eq:glc-condition} and \eqref{eq:llc-condition} are also necessary for such a situation.
Robust adversarial reinforcement learning (RARL) \cite{Pinto2017icml} and its extension \cite{Shioya2018iclr} fall into our framework with $\eta = 1$.
Our theoretical investigation reveals a possible limitation in these two and a need for the introduction of $\eta < 1$.

A local convergence analysis of a simultaneous gradient update \eqref{eq:grad} in \cite{Mescheder2017nips} shows that $(x_t, y_t)$ converges toward $(x^*, y^*)$ if and only if
all the eigenvalues of the matrix $\Big[\begin{smallmatrix}
    I - \eta H_{x,x}^* & - \eta H_{x, y}^* \\
    \eta H_{y, x}^* & I + \eta H_{y,y}^*
  \end{smallmatrix}\Big]$
  live in the unit disk on a complex space.
Again, considering the abovementioned case of $f(x,y) = x^2/2 + \bar\sigma xy - y^2/2$, 
the absolute values of the eigenvalues in a complex space of the above matrix are all $((1-\eta)^2 + \eta^2\bar\sigma^2)^{1/2}$. 
That is, the necessary and sufficient condition $\eta < 2 / (1 + \bar\sigma^2)$ for the simultaneous gradient update \eqref{eq:grad} is equivalent to the condition on $\eta$ for \eqref{eq:algo} with $\bar\epsilon = 0$. 

For the global convergence in \Cref{thm:glc}, we assume that there exists $\delta \in [0, 1)$ satisfying \eqref{eq:delta} for all $(x, y) \in \X\times\Y$.
In contrast, Theorem~1 of \cite{liang2019interaction} shows the global linear convergence of the simultaneous gradient approach with a sufficiently small $\eta$ on a globally and strongly convex--concave $f$.
Our condition $\delta < 1$ is stronger.
We suspect that this is not a potential difference between the simultaneous gradient approach and the zero-order saddle point optimization \eqref{eq:algo} and that it is not necessary for the global convergence itself.
The looseness of the result will be demonstrated in \Cref{sec:exp}.

\section{Learning Rate Adaptation}\label{sec:lra}

Often, $\bar\sigma$ is unknown, especially when zero-order optimization is required.
Here we propose a heuristic approach to adapt the learning rate $\eta$.
We implemented zero- and first-order saddle point optimization algorithms using a randomized hill-climbing algorithm, namely, the (1+1)-ES with 1/5 success rule \cite{Devroye:72,Schumer:Steiglitz:68,rechenberg:1973}, and a sequential least squares programming (SLSQP) subroutine \cite{kraft1988software}. Our implementation is publicly available\footnote{\url{https://gist.github.com/youheiakimoto/15212dbf46dc546af20af38b0b48ff17}}.

\subsection{Adaptation Mechanism}



The suboptimality error $G(x_t, y_t)$ is approximated by using the oracle outputs $x' = M^{A_x}_{\epsilon_x,r_x}(f(\cdot, y_t), x_t)$ and $y' = M^{A_y}_{\epsilon_y,r_y}(-f(x_t, \cdot), y_t)$ as $G(x_t, y_t) \approx f(x_t, y') - f(x', y_t) =: F(x_t, y_t)$.
More precisely, we have
$(1 - 2\bar\epsilon)G(x_t,y_t) \leq F(x_t,y_t) \leq G(x_t,y_t)$. 
This implies that if there exists $\tilde\gamma < 1$ such that $F(x_t, y_t) \leq \tilde\gamma^{2t} F(x_0, y_0)$, we can guarantee that
$G(x_t, y_t) \leq (1 - 2\bar\epsilon)^{-1} \tilde{\gamma}^{2t} G(x_0, y_0)$.
We note also that $G \approx \tilde{G}$ in a neighborhood of $(x^*, y^*)$ for $\mathcal{C}^2$ functions.

\begin{algorithm}[t]\small
  \caption{Learning Rate Adaptation}\label{alg:lr}
  \begin{algorithmic}[1]
    \Require $x\in\X$, $y\in\Y$, $a_\eta > 0$, $b_\eta \geq 0$, $c_\eta > 1$
    \State $\eta \leftarrow 1$, $\tilde\gamma \leftarrow 0$
    \For{$t = 1, \cdots, T$}
    \State $x_t \leftarrow x$, $y_t \leftarrow y$\label{l:keep}
    \State $\eta_c \leftarrow \{ \min(\eta \cdot c_\eta, 1), \eta, \eta / c_\eta \}$ w.p.\ $1/3$ for each
    \State $N_\text{step} \leftarrow \lfloor b_\eta +  a_\eta / \eta_c \rfloor$
    \For{$s = 1, \cdots, N_\text{step}$}\label{l:xcyc}
    \State $\hat{x} \leftarrow M^{A_x}_{\epsilon_x,r_x}(f(\cdot, y), x)$\label{l:xhat}
    \State $\hat{y} \leftarrow M^{A_y}_{\epsilon_y,r_y}(-f(x, \cdot), y)$\label{l:yhat}
    \State $F_{s} \leftarrow  f(x, \hat{y}) - f(\hat{x}, y)$\label{l:gap}
    \State $x \leftarrow x + \eta_c (\hat{x} - x)$\label{l:xup}
    \State $y \leftarrow y + \eta_c (\hat{y} - y)$\label{l:yup}
    \State \textbf{break} \textbf{if} $s \geq b_\eta$ and $F_{s} > \cdots > F_{s-b_\eta+1}$
    \EndFor
    \State $\tilde\gamma_c,\ \sigma_{\tilde\gamma_c} \leftarrow  \textsc{slope}(\log(F_1), \dots, \log(F_{s}))$\label{l:slope}
    \If{$\tilde\gamma \geq 0$ and $\tilde\gamma_c \geq 0$}
    \State $\eta \leftarrow \eta / c_\eta^3$\label{l:etadec}
    \ElsIf{$\tilde\gamma_c \leq \tilde\gamma$ or $\eta = \eta_c$}
    \State $\eta \leftarrow \eta_c$, $\tilde\gamma \leftarrow \tilde\gamma_c$\label{l:etaup}
    \EndIf
    \State $x \leftarrow x_t$, $y \leftarrow y_t$ \textbf{if} $\tilde\gamma - 2 \sigma_{\tilde\gamma_c} > 0$\label{l:revert}
    \EndFor
  \end{algorithmic}
\end{algorithm}

From the above observation, we propose to adapt the learning rate $\eta > 0$ by \Cref{alg:lr}. 
The object is to find the value of $\eta$ that minimizes the convergence rate of $F(x_t, y_t)$.
We perform a random search for the optimal $\eta$. 
The learning rate $\eta_c$ to be tested is one of $\eta$, $\eta \cdot c_\eta$, or $\eta / c_\eta$ chosen at random with a probability of $1/3$ each.
The logarithm of the convergence rate, $\tilde\gamma_c$, is estimated by running the algorithm \eqref{eq:algo} with $\eta_c$ for $N_\text{step}$ iterations (Lines~\ref{l:xcyc} to \ref{l:slope}).
Based on the estimated log convergence rate $\tilde\gamma_c$, we search for the best $\eta$.

To estimate the convergence rate, we run $N_\text{step} \propto 1/\eta$ iterations.
The rationale behind the choice of $N_\text{step} \propto 1/\eta$ is as follows.
From \Cref{thm:glc}, we know that the logarithm of the convergence rate, $\log(\gamma)$ in \eqref{eq:glc-conclusion}, is $- \mathcal{O}(\eta)$. 
On the other hand, because of the approximation error, i.e., $G(x_t, y_t) \leq (1 - 2\bar\epsilon)^{-1} \tilde{\gamma}^{2t} G(x_0, y_0)$,
the logarithm of the convergence rate of $F$ estimated over $N_\text{step}$ iterations can deviate from that of $G$ by $\pm \abs{\log(1 - 2 \bar\epsilon)} / N_\text{step}$.
To alleviate the effect of the approximation error, we need $\epsilon / N_\text{step} \in \mathcal{O}(\eta)$.
In words, $N_\text{step}$ needs to be no less than proportional to $1/\eta$.

\Cref{alg:lr} summarizes the overall framework of the proposed approach with the adaptation mechanism.%
\footnote{In practice, it is advised to try out a fresh $y$ after every $M^{A_y}_{\epsilon_y,r_y}(-f(x_t, \cdot), y_t)$ calls if the minimization oracle is for a local minimization but it is desired to avoid (rather) suboptimal $y$.
  After Line~\ref{l:yhat}, one may add the following lines:
  \begin{algorithmic}[1]
    \State $\tilde{y} \sim p_{y}$
    \If{ $f(x, \tilde{y}) > f(x, \hat{y})$ }
    \State $\hat{y} \leftarrow \tilde{y}$
    \EndIf
  \end{algorithmic}
  A possible choice of $p_y$ is the distribution from which the initial $y$ is drawn.
  If the search space of $y$ is bounded, the uniform distribution over the search domain of $y$ is a candidate for $p_y$.
  One can also perform analogous steps for $x$.
  We didn't implement it in our experiments as the objective of the experiments are to understand the proposed approach and the adaptation mechanism.
}
We introduce three hyperparameters $a_\eta$, $b_\eta$, and $c_\eta$ to adapt one parameter $\eta$.
Arguably, tuning them is easier than tuning $\eta$ itself.
A greater $a_\eta$ and $b_\eta$ will lead to a more accurate estimation of $\tilde\gamma_c$, while spending more $f$-calls and requiring a longer adaptation time for $\eta$.
The granularity of the randomized line search for $\eta$ is controlled by $c_\eta$.
A smaller $c_\eta$ will lead to a smoother change in $\eta$, while requiring a longer adaptation time.

\begin{algorithm}[t]\small
  \caption{(1+1)-ES}\label{alg:es}
  \begin{algorithmic}[1]
    \Require $h: \R^\ell \to \R$, $z \in \R^\ell$, $\sigma \in (0, \sigma_{\max}]$, $h_{z} = h(z)$, $\tau_\text{es} \in \mathbb{N}$
    \State $c = e^{1/\sqrt{2 \ell}}$
    \State $n_\text{succ} = 0$
    \While{$n_\text{succ} < \tau_\text{es} \cdot \ell$}
    \State\label{alg:aes:esstart} $z' \leftarrow z + \sigma \mathcal{N}(0, I)$
    \State $h_{z'} \leftarrow h(z')$
    \If{$h_{z'} \leq h_{z}$} 
    \State $\sigma \leftarrow \min(\sigma \cdot c, \sigma_{\max})$, $z \leftarrow z'$, $h_{z} \leftarrow h_{z'}$
    \State $n_\text{succ} \leftarrow n_\text{succ} + 1$
    \Else
    \State $\sigma \leftarrow \sigma \cdot c^{-1/4}$
    \EndIf\label{alg:aes:esend}
    \EndWhile
    \State \Return $z$, $h_z$, $\sigma$
  \end{algorithmic}
\end{algorithm}

\begin{figure*}[t]
  \centering
  \begin{subfigure}{0.33\hsize}%
    \centering%
    \includegraphics[width=\hsize]{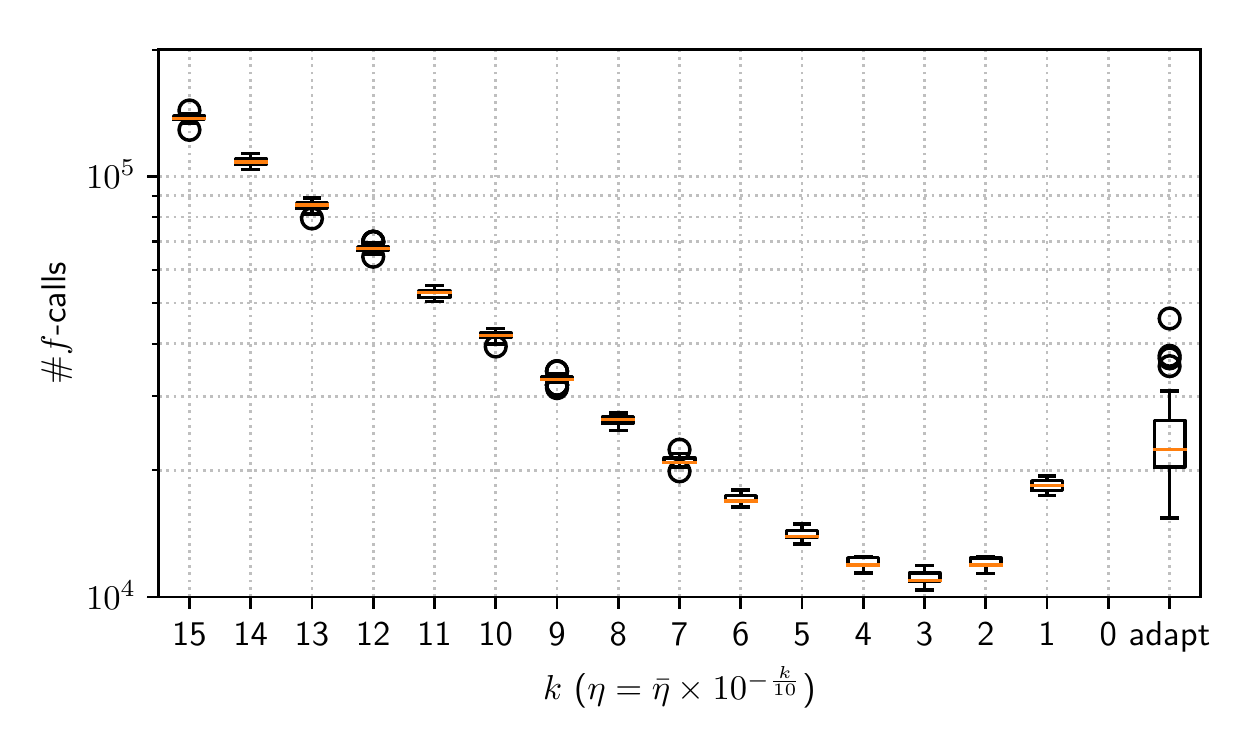}%
    \caption{AES ($n = 10$, $b = 1$)}%
  \end{subfigure}%
  \begin{subfigure}{0.33\hsize}%
    \centering%
    \includegraphics[width=\hsize]{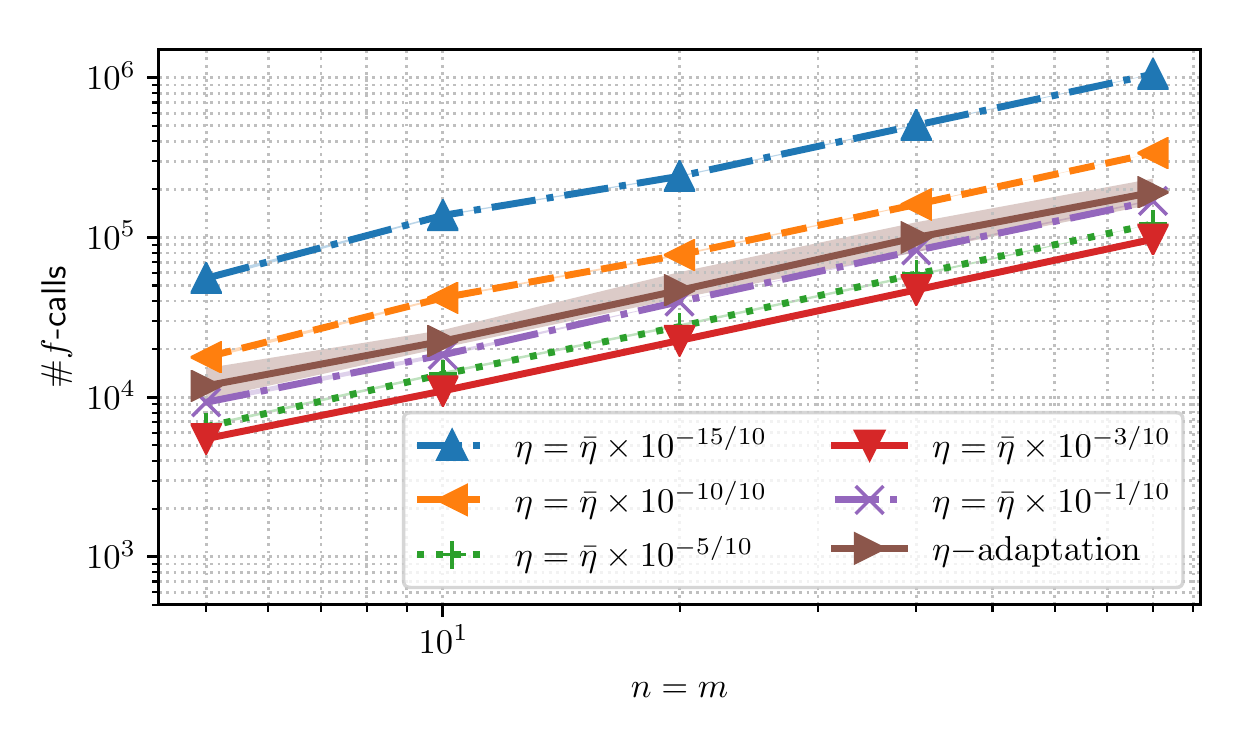}%
    \caption{AES (varying $n$, $b = 1$)}%
  \end{subfigure}
  \begin{subfigure}{0.33\hsize}%
    \centering%
    \includegraphics[width=\hsize]{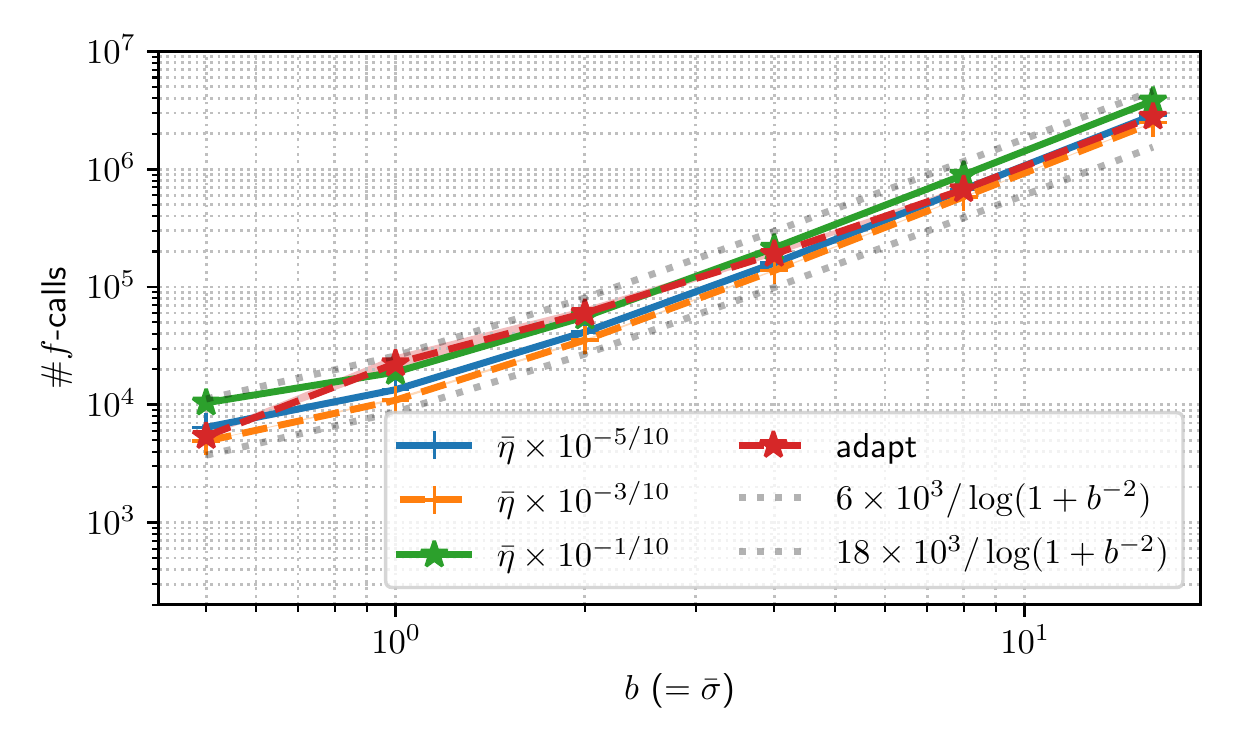}%
    \caption{AES ($n =10$, varying $b$)}%
  \end{subfigure}%
  \\
    \begin{subfigure}{0.33\hsize}%
    \centering%
    \includegraphics[width=\hsize]{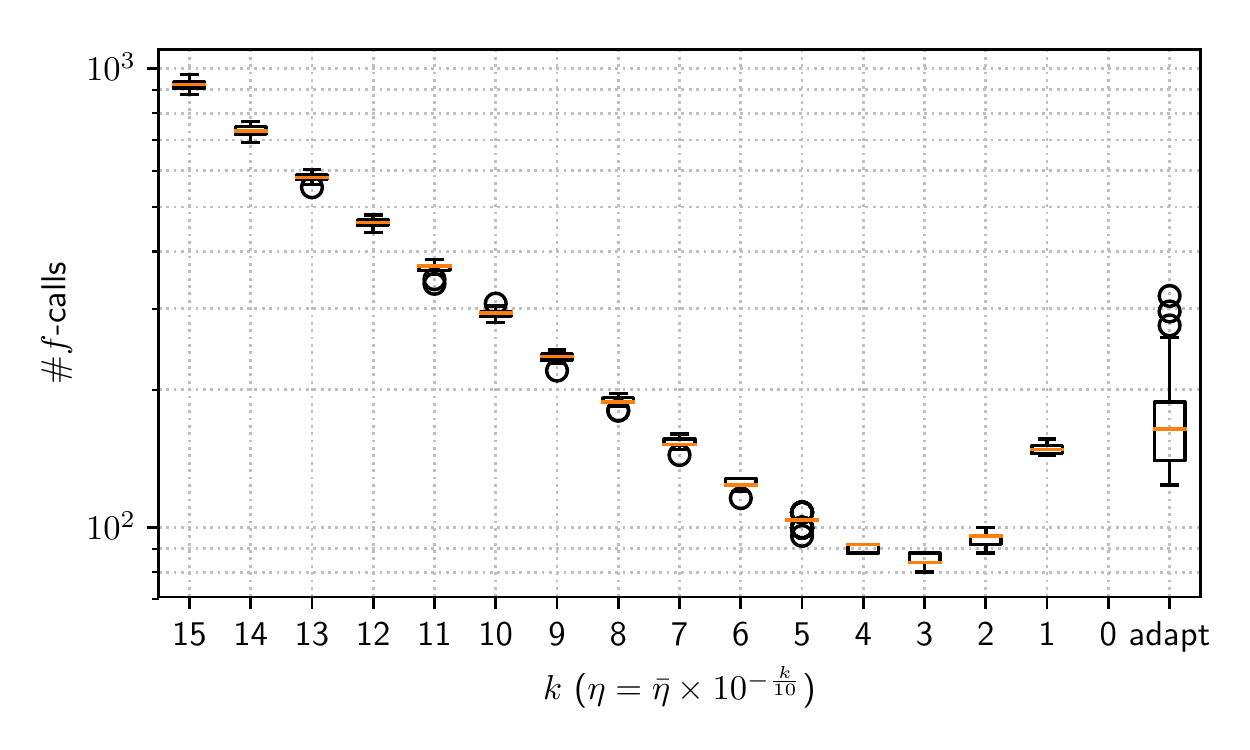}%
    \caption{ASLSQP ($n = 10$, $b = 1$)}%
  \end{subfigure}%
  \begin{subfigure}{0.33\hsize}%
    \centering%
    \includegraphics[width=\hsize]{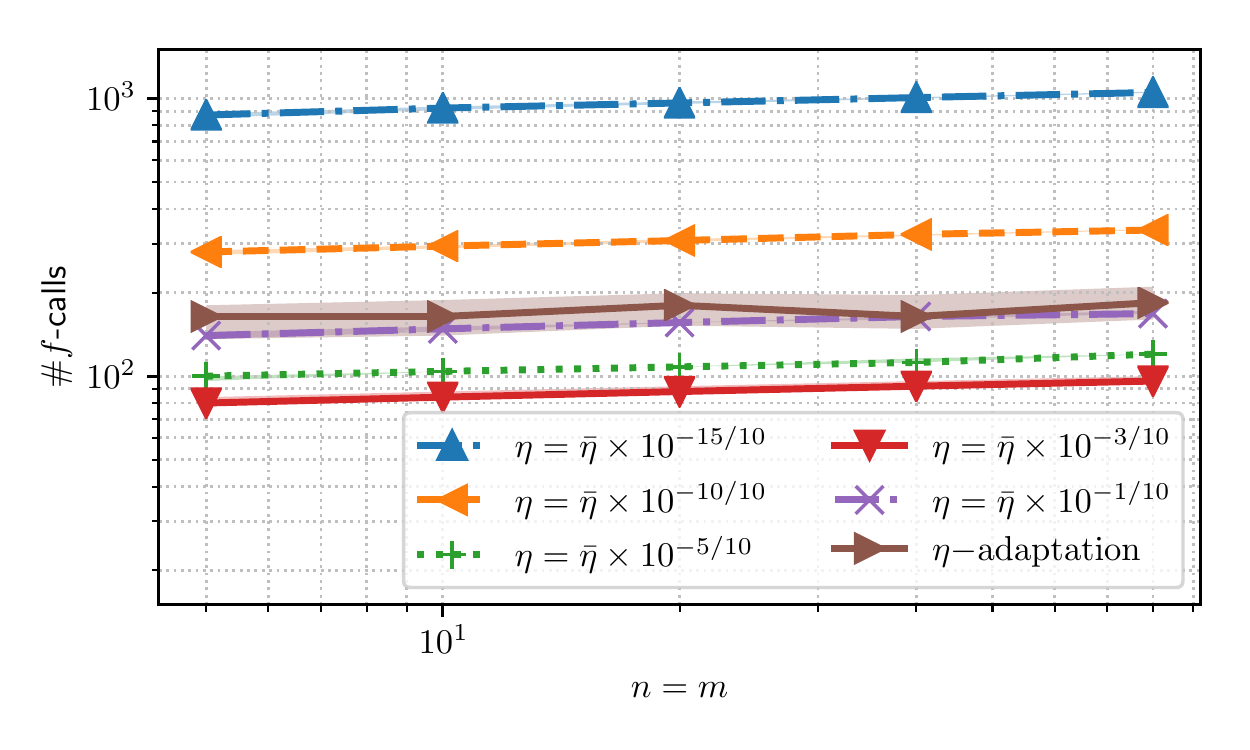}%
    \caption{ASLSQP (varying $n$, $b = 1$)}%
  \end{subfigure}%
  \begin{subfigure}{0.33\hsize}%
    \centering%
    \includegraphics[width=\hsize]{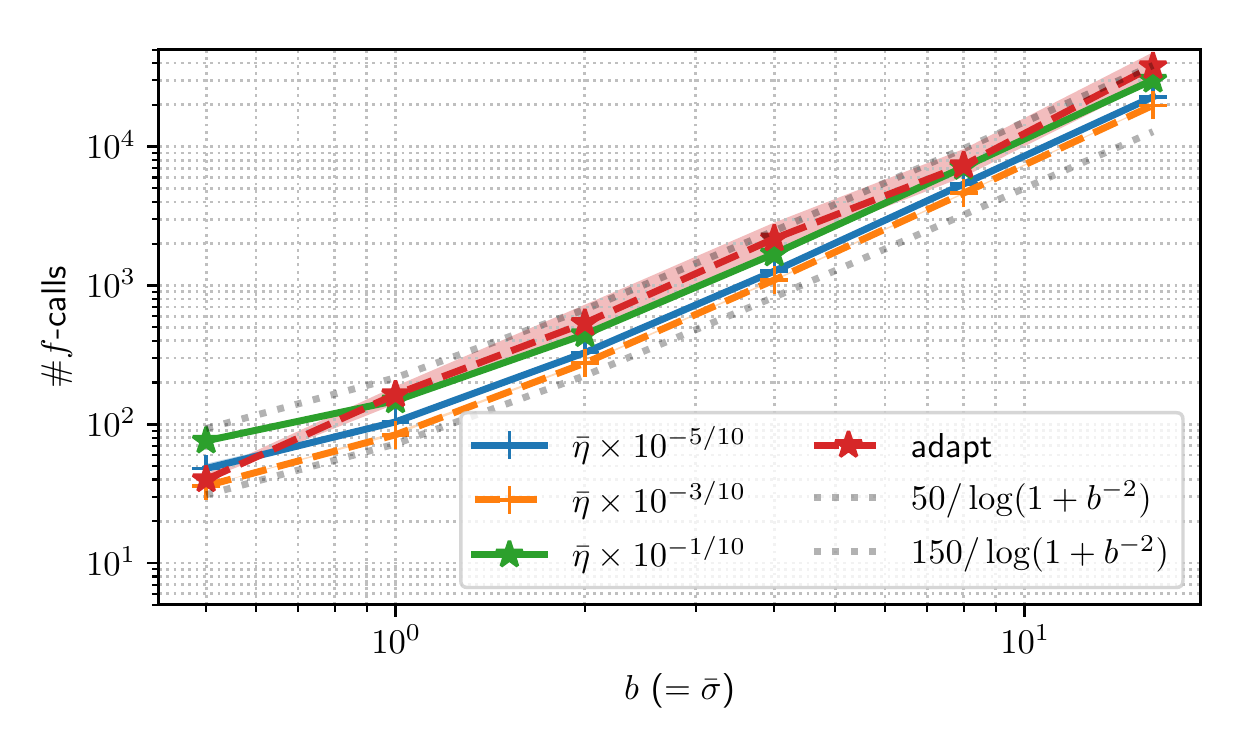}%
    \caption{ASLSQP ($n = 10$, varying $b$)}%
  \end{subfigure}%
  \caption{The number of $f$-calls until $G(x_t, y_t) \leq 10^{-5}$ is obtained on $f_1$. (a--c) AES. (d--f) ASLSQP.
    (a, d) Boxplots of the results with and without $\eta$ adaptation.
    Missing data imply that no run has succeeded in achieving $G(x_t, y_t) \leq 10^{-5}$.
    (b, e) Median (solid line) and inter-quartile range (band) for dimensions $n=m$. 
    (c, f) Median (solid line) and inter-quartile range (band) for interaction term $b$. 
    }
\label{fig:ex1}%
\end{figure*}

\subsection{Adversarial Evolution Strategy}

We approximate $M_{\epsilon,r}^A$ with the so-called (1+1)-ES with a $1/5$-success rule.
The (1+1)-ES is a randomized hill-climbing algorithm with an adaptive step-size mechanism.
We adopt the simplified update of the step-size proposed by \cite{kern:2004}.
\Cref{alg:es} summarizes the (1+1)-ES on $h: \Z \to \R$.
In our case, $h$ corresponds to $f(\cdot, y_t)$ for $x$ updates and $- f(x_t, \cdot)$ for $y$ updates.
It generates a candidate solution $z'$ from a normal distribution centered at $z$ with standard deviation $\sigma$.
If the candidate solution has a function value smaller than or equal to that of the current solution, we replace the current solution with the candidate and call each such iteration successful.
The step size is increased by multiplying $c$ if the iteration is successful.
Otherwise, it is decreased by multiplying $c^{-1/4}$.
Such a step-size adaptation mechanism is called a $1/5$-success rule because it maintains the step size to keep the probability of success to $1/5$.

The runtime time analysis of the (1+1)-ES algorithm reveals that the expected number $\E[T]$ of iterations $T$ to achieve $\norm{z - z^*} \leq \epsilon R_0$ from the initial condition $R_0 = \norm{z - z^*}$ is $\Theta(\log(1/\epsilon))$ on strongly convex functions with Lipschitz continuous gradients and their strictly increasing transformations \cite{morinaga2019generalized}.
Moreover, its scaling w.r.t.~dimension $\ell$ is proven to be $\Theta(\ell)$ on convex quadratic functions \cite{morinaga2021drift}.
Therefore, we expect that \Cref{alg:es} approximates $M_{\epsilon,r}^A$ with $\epsilon \in \exp(-\Theta(T / \ell))$.

We propose a zero-order saddle point optimization algorithm using the (1+1)-ES, called \emph{Adversarial Evolution Strategy (AES)}.
AES replaces $M^{A_x}_{\epsilon_x,r_x}(f(\cdot, y_i), x_i)$ and $M^{A_y}_{\epsilon_y,r_y}(-f(x_i, \cdot), y_i)$ in \Cref{alg:lr} with \Cref{alg:es}.
In \Cref{alg:lr}, in addition to $x$ and $y$, AES maintains two step-size parameters $\sigma_x$ for $x$ update and $\sigma_y$ for $y$ update.
On Lines~\ref{l:xhat} and \ref{l:yhat}, we input $\sigma_{x}$ and $\sigma_{y}$ to \Cref{alg:es} to obtain the updated values.
On Lines~\ref{l:keep} and \ref{l:revert}, we maintain $\sigma_x$ and $\sigma_y$ as well as $x$ and $y$.

Importantly, the number $T$ of iterations of \Cref{alg:es} is not fixed for each oracle call.
Instead, we count the number of successful iterations, where a better or equally accurate objective function value is generated.
The reason for this design choice is as follows.
As mentioned above, we expect $\epsilon \in \exp(-\Theta(T / \ell))$.
However, if, for example, $\sigma$ is too large to produce a successful candidate solution\footnote{It is reported in \cite{morinaga2021drift} that $\sigma$ needs to be proportional to the norm of the gradient to produce a successful candidate solution if the objective function is convex quadratic and that the (1+1)-ES with 1/5 success rule maintains $\sigma$ to be proportional to the norm of the gradient.}, \Cref{alg:es} has little chance to produce any improvement until $\sigma$ is well adapted.
Because the correct $\sigma$ approaches zero as the candidate solution approaches the optimum, it requires more iterations to adapt $\sigma$ if it is initialized to a fixed value for each oracle call. 
Then, the output of \Cref{alg:es} with a fixed number of iterations will not satisfy the oracle condition \eqref{eq:oracle} for some $\epsilon$.
To reduce the effect of an excessively large initial $\sigma$, we count the number of successful iterations rather than total iterations.

To further relax the effect of an unsuitable initial $\sigma$ value, AES shares $\sigma$ for successive oracle calls.
With this parameter being shared between oracle calls, we expect the following.
Compared to the case where the initial $\sigma$ is too large, we expect this to reduce the number of iterations (objective function calls) to satisfy the stopping condition because the adaptation time for $\sigma$ is reduced.
Compared to the case where the initial $\sigma$ is too small, where the iteration can be successful with high probability despite the improvement being rather small, we expect to avoid excessively early stopping of \Cref{alg:es}.\footnote{Our preliminary experiments  revealed that the $\sigma$ sharing effect is statistically significant; however, it may not be practically important.
  The ratio between the numbers of $f$-calls to reach the same target threshold with AES and AES without $\sigma$ sharing was less than 1.1 on $f_1$ with $n = m = 10$ and $b = 2$.
  See \Cref{sec:exp} for the detailed experimental setting.
  The reason is simply that $\sigma$ adaptation in the (1+1)-ES is rather quick.
  From \cite{akimoto2018drift}, we know that the adaptation to forget an inaccurate initial step size is proportional to $\log(\sigma_{0} / \sigma^*)$, where $\sigma_0$ is the initial step size and $\sigma^*$ denotes the optimal step size, which is proportional to $\norm{x_0 - x^*}$ on the spherical function with optimum at $x^*$.
}

\subsection{Adversarial SLSQP}

To demonstrate the applicability of the proposed $\eta$ adaptation mechanism, we also implement it with SLSQP as $M_{\epsilon,r}^A$.
This is a sequential quadratic programming approach. 
We allow access to the gradient of $f$ to the SLSQP procedure and limit the maximum number of iterations to $\tau_\mathrm{slsqp}$.
In contrast to AES, we share nothing but the solutions $x$ and $y$ between different oracle calls.
We call this first-order approach \emph{Adversarial SLSQP (ASLSQP)}.

\section{Experiments}\label{sec:exp}

We conducted numerical experiments with AES and ASLSQP with and without $\eta$ adaptation.
Our objective was to confirm the tightness of the bounds obtained in \Cref{thm:glc,thm:llc} and the applicability of the proposed $\eta$ adaptation mechanism.
AES and ASLSQP without $\eta$ adaptation repeat Lines~\ref{l:xhat} to \ref{l:yup} with a fixed learning rate.
We set $a_\eta = 1$, $b_\eta = 5$, $c_\eta = 1.1$, and $\tau_\text{es} = \tau_\text{slsqp}= 5$ for the following experiments, unless otherwise specified.

\subsection{Ex 1: Convex--concave Quadratic Case}\label{sec:ex1}

In this experiment, we confirm that (A) $\bar\eta$ in \eqref{eq:glc-condition} is a tight upper bound; (B) $\bar\eta^*$ in \eqref{eq:optimal} is a good estimate of the optimal $\eta$; (C) $\bar\eta$ in \eqref{eq:glc-condition} and $\eta^*$ and $\bar\gamma^*$ in \eqref{eq:optimal} accurately reflect the dependency of $\bar\sigma$; and (D) the proposed $\eta$ adaptation demonstrates a reasonable performance compared to the best fixed $\eta$.

We consider the following convex--concave quadratic function $f: \R^m \times \R^n \to \R$:
\begin{equation*}
  f_1(x, y) = \norm{x}^2/2 + b x^\T y - \norm{y}^2/2 \enspace.
\end{equation*}
For this function, $\bar\sigma = b$ and $\delta = 0$ in \Cref{thm:glc}.
The upper bound on $\eta$ for linear convergence is $\bar\eta = 2 / (1 + b^2)$.
Moreover, $\bar\gamma^* = \abs{b} / (1 + b^2)^{1/2}$ and $\eta^* = 1/(1 + b^2)$ if $\bar\epsilon = 0$\footnote{The reason that we approximated $\bar\gamma^*$ and $\eta^*$ with $\epsilon = 0$ in the experiments is because the rigorous formula \eqref{eq:optimal} is less understandable. In the experiments, we observed that the oracle condition \eqref{eq:oracle} is satisfied with $\epsilon < 5 \times 10^{-4}$ for 50\% of iterative checks and $\epsilon < 5 \times 10^{-3}$ for 99\% of iterative checks for AES with $\tau = 5$ on $n = 20$, and even smaller for ASLSQP. Such small values \eqref{eq:optimal} can be well approximated by those computed with $\bar\epsilon = 0$.}.
We measured the performance of the algorithms by the number of $f$-calls until they reached $G(x_t, y_t) \leq 10^{-5}$.
For each setting, we ran $50$ independent trials with different initial search points generated from $\mathcal{N}(0, I)$.
The step-size $\sigma_x$ and $\sigma_y$ for AES are initialized to $2$ and their maximal values, denoted by $\sigma_{\max}$ in \Cref{alg:es}, are also set to $2$.

To confirm (A--C), we ran algorithms without $\eta$ adaptation with $\eta = \bar\eta \times 10^{-k/10}$ for $k = 0,1,\dots,15$ on $f_1$ with (i) $n = m = 10$ and $b = 0.5$, $1$, $2$, $4$, $8$, $16$ and (ii) $n = m = 5$, $10$, $20$, $40$, $80$ and $b = 1$. To confirm (D), we ran algorithms with $\eta$ adaptation on the same problems. The results are summarized in \Cref{fig:ex1}.

(A) No run succeeded with $\eta = \bar\eta$ for any $n$ and $b$, as demonstrated in \Cref{fig:ex1}-(a, d) for the case of $n = m = 10$ and $b = 1$, whereas all runs succeeded with $\eta = \bar\eta \times 10^{-1/10}$. This demonstrates the tightness of $\bar\eta$.

(B) For all tested cases, $\eta = \bar\eta \times 10^{-3/10} \approx \eta^*$ showed the best performance, as demonstrated in \Cref{fig:ex1}-(a--f). That is, $\eta^*$ in \eqref{eq:optimal} is a good estimate of the best $\eta$.

(C) In \Cref{fig:ex1}-(c, f), we observe that the best $\eta$ was approximately $\bar\eta\times 10^{-3/10} \approx \eta^*$ for all $b$ values. The number of $f$-calls was proportional to $1 / \log(1 + b^{-2}) \in \mathcal{O}(1/\log(1/\bar\gamma^*))$, implying that our expectation of the dependency of $\bar\sigma = b$ on $\gamma^*$ is accurate. Moreover, as mentioned above, we observe that all runs succeeded with $\eta = \bar\eta \times 10^{-1/10}$, whereas no run was solved with $\eta = \bar\eta$. We conclude that $\bar\eta$, $\eta^*$, and $\bar\gamma^*$ accurately reflect their dependent value $\bar\sigma$ as expected.



(D) \Cref{fig:ex1}-(a--f) shows that AES and ASLSQP with $\eta$ adaptation solved the problems with up to three times more $\#f$-calls than AES and ASLSQP with the best fixed $\eta = \bar\eta \times 10^{-3/10}$ for upper quartile cases.
We emphasize that for AES without $\eta$ adaptation to achieve performance comparable to AES with $\eta$ adaptation, one needs to find $\eta \in [\bar\eta/10, \bar\eta)$. 
It is rather difficult to find such intervals within a few trials when $\bar\sigma$ is unknown.

\Cref{thm:glc,thm:llc} indicate no dependency of $n$ and $m$ on $\gamma$.
This implies that the number of oracle calls does not scale with $n$ and $m$, whereas
the number of $f$-calls may depend on them.
In ASLSQP (\Cref{fig:ex1}-(e)), $\#f$-calls is nearly constant for varying dimensions as its runtime does not scale with dimension. (note that ASLSQP also requires a similar amount of $\#\nabla f$-calls), whereas AES (\Cref{fig:ex1}-(b)) requires $\#f$-calls proportional to the dimension, as its runtime is proportional.

\subsection{Ex 2: Convex--concave Case with $\delta > 1$}

Next, we demonstrate that $\delta < 1$ in \Cref{thm:glc} is not necessary for convergence.
We consider the following function $f: \R^m \times \R^n \to \R$ with $m = n = 10$:
  \begin{equation*}
    f_2(x, y) = f_1(x, y) - \exp(- \norm{x}^2/2) + \exp(- \norm{y}^2/2)
  \end{equation*}
  On such a function, $H_{x,x} = (1 + (1 - x^2)\exp(- x^2/2)) I$, $- H_{y,y} = (1 + (1 - y^2)\exp(- y^2/2)) I$, $H_{x,y} = H_{y,x} = b I$ and $H_{x,x}^* = - H_{y,y}^* = 2 I$. Then, $G_{x,x}^* = G_{y,y}^* = ((4 + b^2) / 2) I$.
  We have $\delta > b / 2$. 
  Therefore, $b = 10$ violates the condition for \Cref{thm:glc}.
  For $b = 10$, ignoring the effect of $\delta$ (i.e., considering local convergence), we have $\bar\eta = 4 / (4 + b^2)$ and  $\eta^* = 2/(4 + b^2)$ if $\bar\epsilon = 0$.
  We conducted experiments with the same setting as \Cref{sec:ex1}.
  We measure the progress by the approximate suboptimality error $\tilde{G}$ around the saddle point.

\begin{figure}[t]
  \centering
  \includegraphics[width=0.9\hsize]{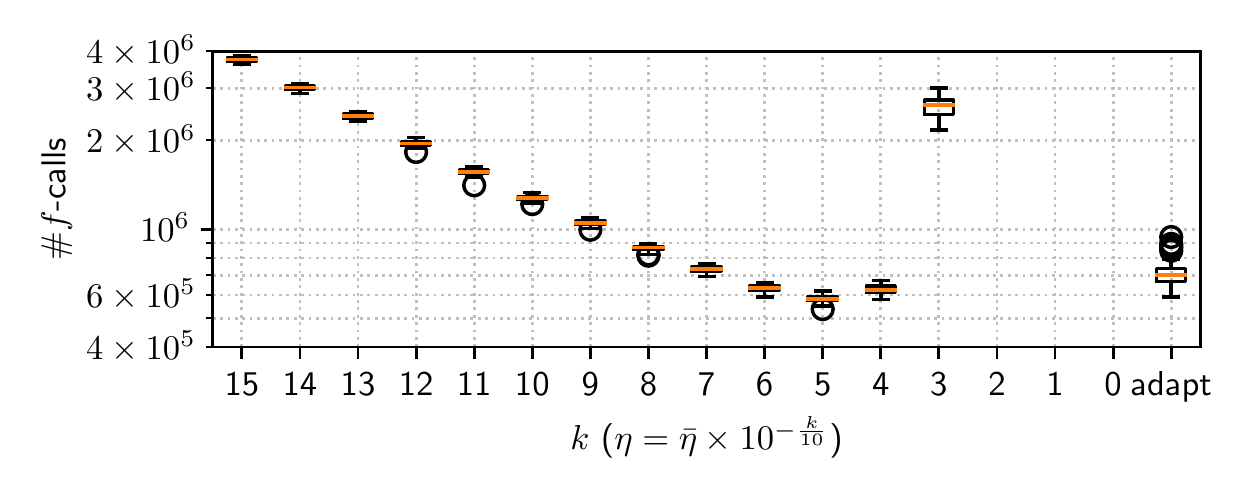}\\%
  \includegraphics[width=0.9\hsize]{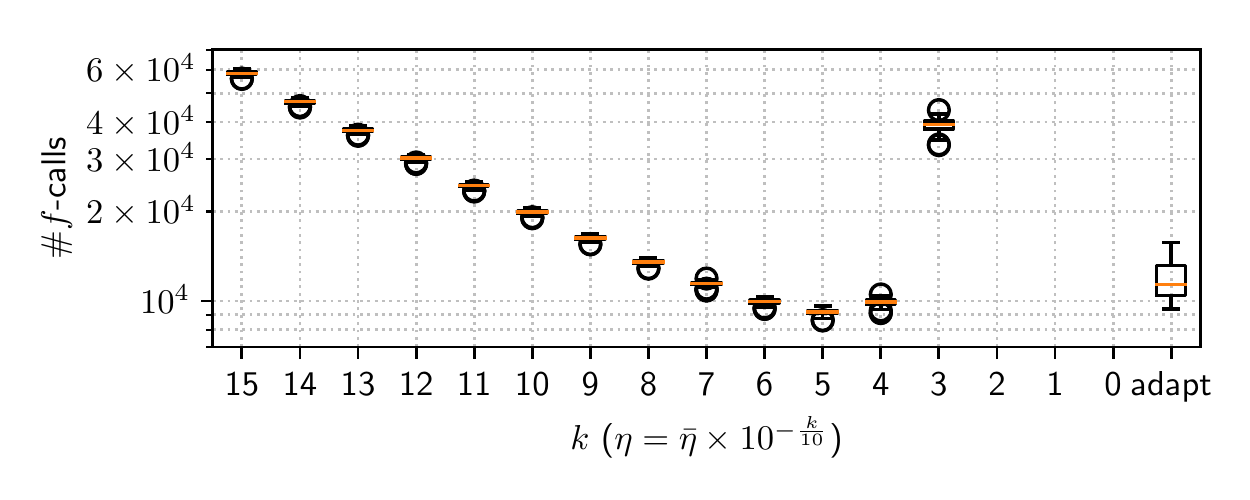}%
\caption{The number of $f$-calls until $\tilde{G}(x_t, y_t) \leq 10^{-5}$ is obtained with AES (top) and ASLSQP (bottom) with $\eta = \bar\eta \times 10^{-k/10}$ for $k = 0, \dots, 15$ and $\eta$ adaptation on $f_2$ with $b = 10$ and $n = m = 10$.}
\label{fig:ex2}%
\end{figure}
  
\Cref{fig:ex2} shows the boxplot of the number of $f$-calls until $\tilde{G}(x_t, y_t) \leq 10^{-5}$ is reached.
As suggested in \Cref{sec:discussion}, both algorithms were able to locate a nearly optimal saddle point with a sufficiently small $\eta$.
This reveals room for improvement in \Cref{thm:glc}.
We conjecture that the upper bound on $\eta$ to guarantee linear convergence will be similar to those obtained in \cite{liang2019interaction} for the simultaneous gradient update. 

\subsection{Ex 3: NonConvex--concave Case}

\begin{figure}[t]
  \centering
  \begin{subfigure}{0.5\hsize}%
    \includegraphics[width=\hsize]{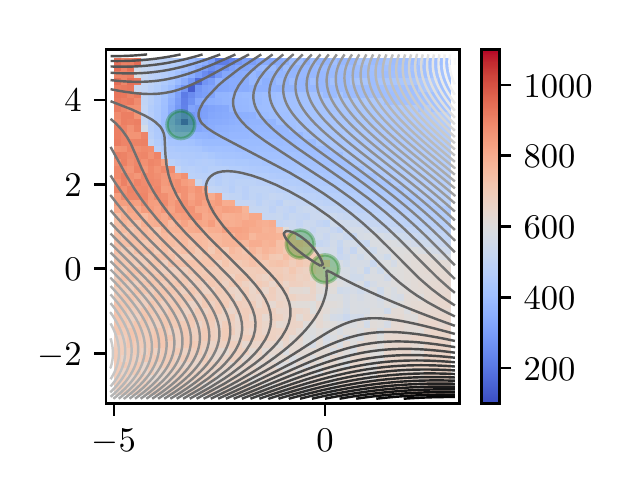}%
    \caption{$\eta = 0.1$, $\tau_\mathrm{slsqp}=1$}
  \end{subfigure}%
  \begin{subfigure}{0.5\hsize}%
    \includegraphics[width=\hsize]{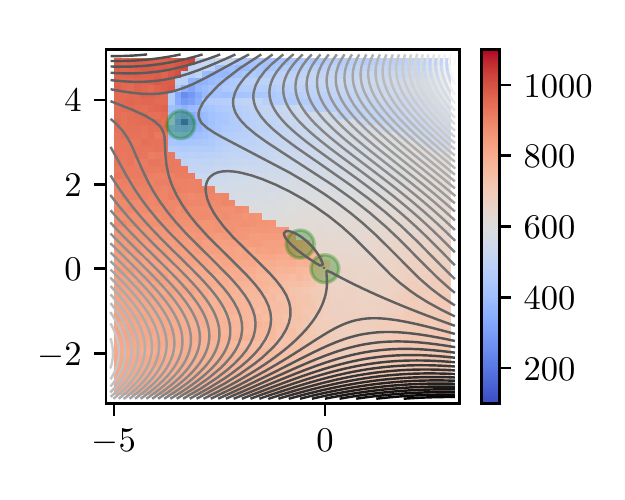}%
    \caption{$\eta = 0.1$, $\tau_\mathrm{slsqp}=5$}    
  \end{subfigure}%
  \caption{$\#f$-calls until $\tilde{G}(x_t, y_t) \leq 10^{-5}$ is obtained from each initial point.
  The three circles are $z_1$, $z_2$ and $z_3$.}
\label{fig:ex3}%
\end{figure}

A function with nonoptimal critical points
\begin{equation*}
  f_3(x, y) = 2 x^2 + 4xy + y^2 + (4/3)y^3 - (1/4)y^4 
\end{equation*}
has three critical points $z_0 = (0, 0)$, $z_1 = (-2 - \sqrt{2}, 2 + \sqrt{2})$, and $z_2 = (-2 + \sqrt{2}, 2 - \sqrt{2})$, of which only $z_1$ is the locally optimal saddle point, whereas the others are local minima. 
The simultaneous gradient method \eqref{eq:grad} is reportedly attracted by local minima \cite{pmlr-v89-adolphs19a}.
ASLSQP is also gradient based, but the update step is different from \eqref{eq:grad}, and it performs $\tau_\mathrm{slsqp}$ steps for each oracle call.
A question arises as to whether the same undesired convergence is observed for ASLSQP.
To investigate this question, we run ASLSQP with $\eta = 0.1$ and $\tau_\mathrm{slsqp} = 1$ (single step) and $5$ (multi-step) from $51^2$ different initial points on $[-5, 3]\times[-3, 5]$.
The progress is measured by $\tilde G(x, y)$ around the saddle point $z_1 = (x^*, y^*)$. 
The results are shown in \Cref{fig:ex3}.
All the trials succeeded in locating an approximate saddle point.
From this result, we conjecture that nonsaddle critical points are not attractors of ASLSQP, which is, however, not covered by \Cref{thm:llc}. Future work in this line is required.

\subsection{Ex 4: Limitation of Existing Coevolutionary Approaches}

Typical existing coevolutionary approaches do not fit into our framework \eqref{eq:algo}.
Therefore, the limitations of existing coevolutionary approaches are not revealed in our analysis.
To show the nonconvergent behavior of these approaches, we run \texttt{mmde} \cite{mmde2018}, \texttt{coeva}, and \texttt{reckless} \cite{Al-Dujaili2019lego} on $f_1$ on $m = n = 2, 10$ with $b = 0, 2$. We used the implementation of these three approaches published by the authors of \cite{Al-Dujaili2019lego}. Because their approaches assume a rectangular search space $[0, 1]^m$ and $[0, 1]^n$, the function is shifted such that the optimal saddle point is at $x^* = (0.5, \cdots, 0.5)$ and $y^* = (0.5, \cdots, 0.5)$.
We also run AES with $\eta$-adaptation, where the initial solution is generated uniform-randomly in $[0, 1]^{m} \times [0, 1]^n$, $\sigma_x$ and $\sigma_y$ are initialized to $1/4$, and their maximum values are set to $1$. The box constraint is treated with the mirroring technique, i.e., a solution outside the feasible domain is mapped by applying the transformation $x \mapsto 1 - \abs{\text{mod}( x, 2 ) - 1}$ for each coordinate.

\begin{figure}[t]
  \centering
  \begin{subfigure}{0.5\hsize}%
    \includegraphics[width=\hsize]{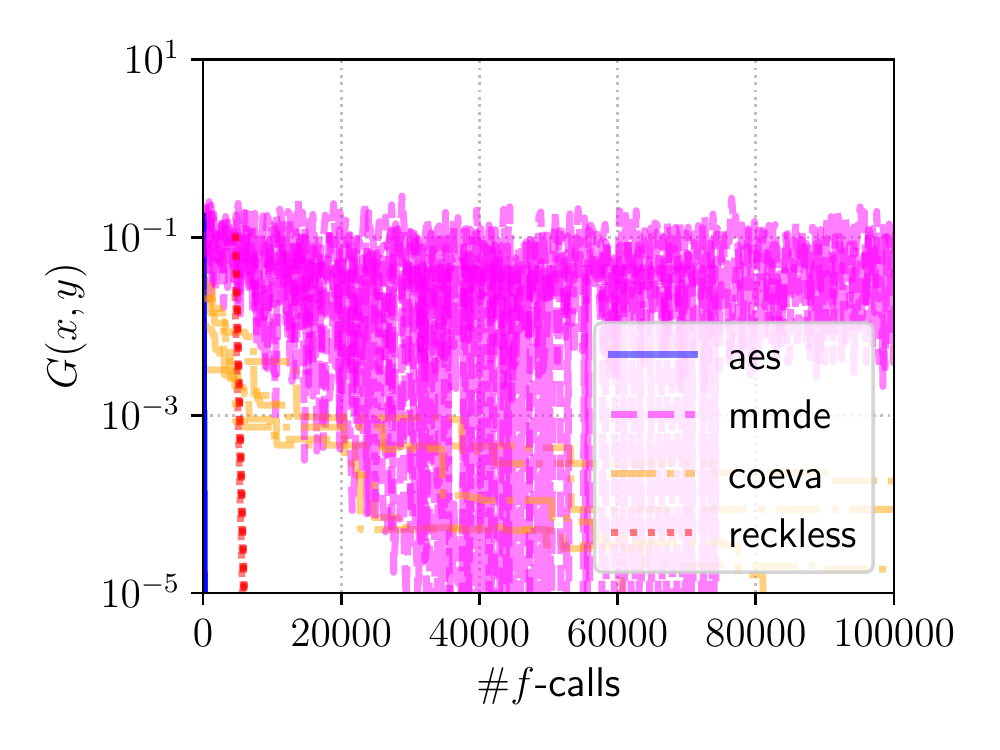}%
    \caption{$b = 0$, $n = m = 2$}
  \end{subfigure}%
  \begin{subfigure}{0.5\hsize}%
    \includegraphics[width=\hsize]{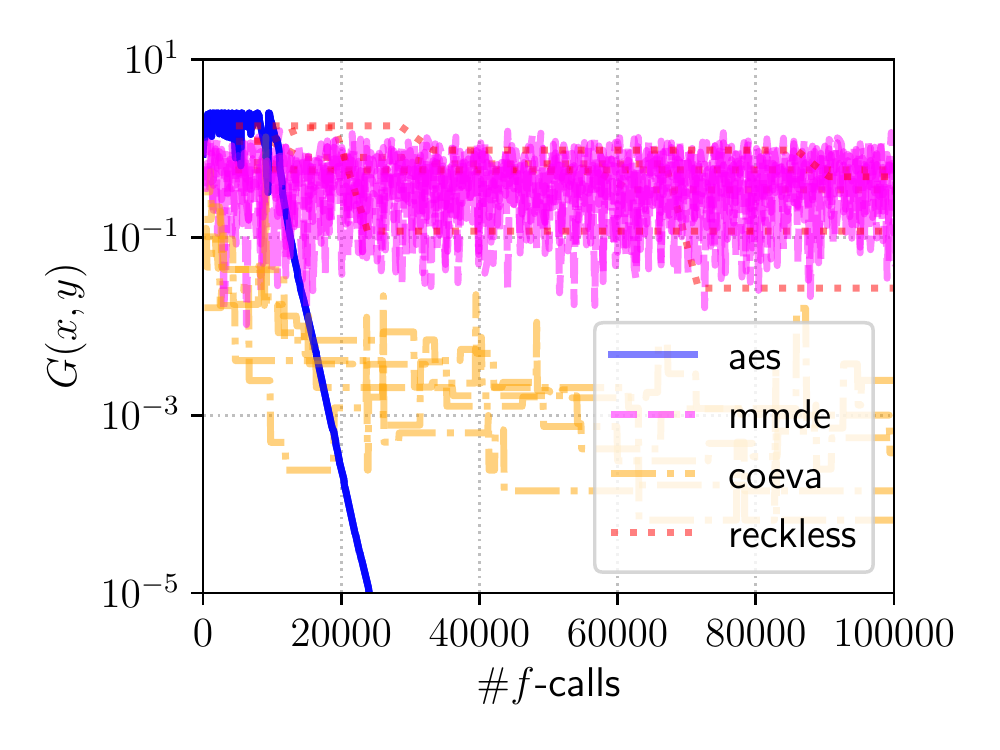}%
    \caption{$b = 2$, $n = m = 2$}
  \end{subfigure}%
  \\%
  \begin{subfigure}{0.5\hsize}%
    \includegraphics[width=\hsize]{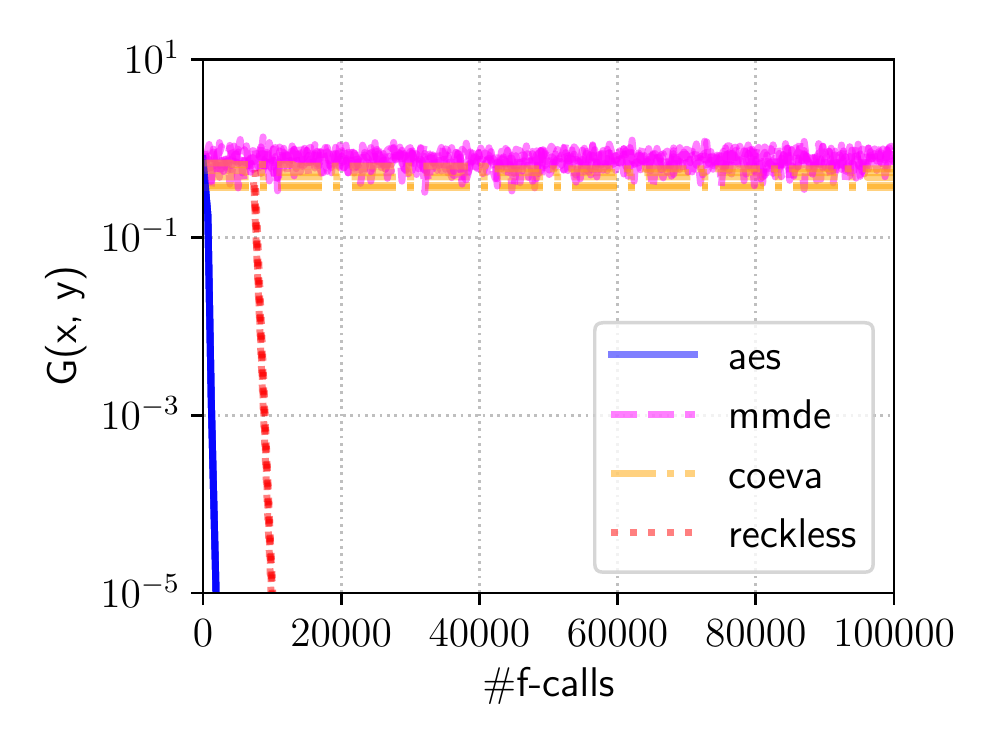}%
    \caption{$b = 0$, $n = m = 10$}
  \end{subfigure}%
  \begin{subfigure}{0.5\hsize}%
    \includegraphics[width=\hsize]{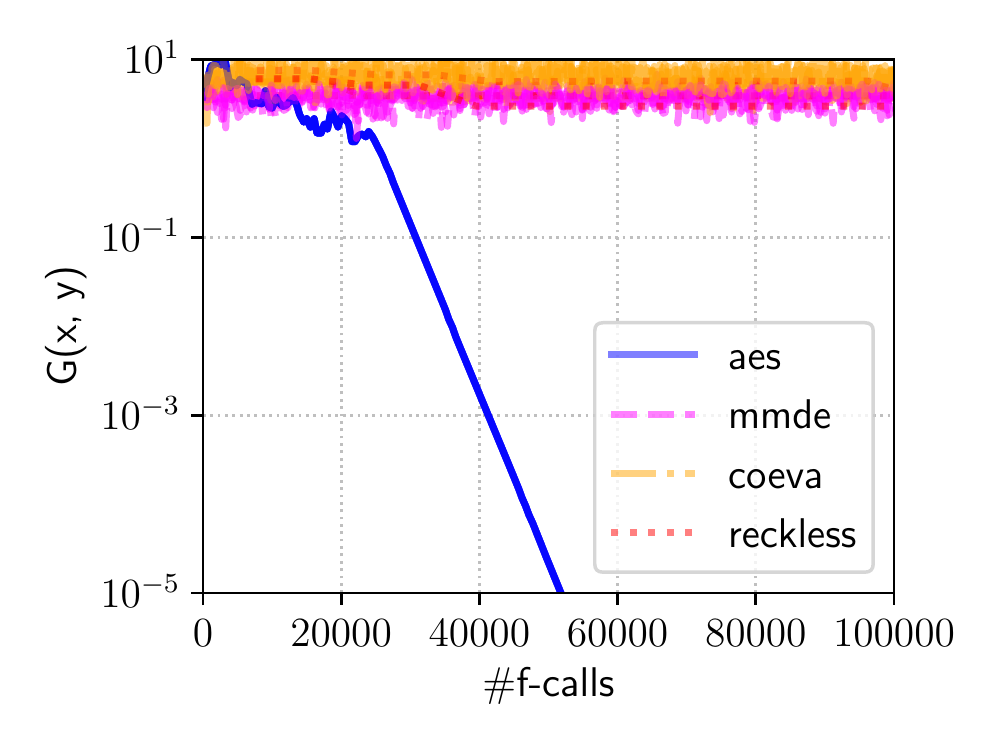}%
    \caption{$b = 2$, $n = m = 10$}
  \end{subfigure}%
  \caption{5 runs of existing coevolutionary approaches on $f_1$}
\label{fig:ex4}%
\end{figure}

The results are shown in \Cref{fig:ex4}.
AES converged linearly towards the global min--max saddle point.
Except for \texttt{reckless} on $b = 0$, the existing approaches, \texttt{reckless}, \texttt{mmde}, and \texttt{coeva}, failed to converge. 
The case $b = 0$ means that there is no interaction between $x$ and $y$, and these variables can be optimized separately. 
Even in such an easy situation, \texttt{mmde} and \texttt{coeva} failed to reach the target quality within $5\times 10^{5}$ $f$-calls. For cases with an interaction term ($b > 0$), no existing algorithms reached the target quality within a given budget, whereas AES reached it with $f$-calls less than $10^{5}$. This demonstrates the difference in convergence behavior between the proposed framework and existing coevolutionary approaches.

\section{Summary and Final Remarks}

A saddle point optimization \eqref{eq:algo} based solely on an approximate minimization oracle is studied.
Its convergence properties are analyzed in \Cref{thm:glc,thm:llc}, where a sufficient condition on the learning rate $\eta$ to exhibit linear (geometric) convergence of an approximate suboptimality error is derived.
One important remark is that simultaneous minimization approaches (\eqref{eq:algo} with $\eta = 1$) do not satisfy the sufficient condition. This may reveal a shortcoming of existing approaches employing simultaneous or alternating updates with $\eta  = 1$ \cite{Pinto2017icml,Shioya2018iclr,Al-Dujaili2019lego}.  
As discussed in \Cref{sec:discussion}, this sufficient condition is also necessary for a convex--concave quadratic function.
Hence, it implies that they do not converge in such a situation and suggests a modification to these algorithms.
A learning rate adaptation heuristic is then proposed and evaluated on test problems.


The generality of our algorithmic framework \eqref{eq:algo} is one of its advantages over analyses on specific algorithms.
The minimization oracle can be zero-order, first-order, etc.
Moreover, it can share internal information during different oracle calls, as it does in AES.
The latter point may be useful for understanding the effect of introducing target networks in deep actor--critic algorithms \cite{lillicrap2015continuous}, which may be regarded as \eqref{eq:algo}, where each oracle call starts from the final solution of the last oracle call.
Further investigation in this direction may be an interesting research topic.

We used the (1+1)-ES algorithm as the minimization oracle in AES for its simplicity to demonstrate the main idea and its theoretical aspects. 
For practical use, it is generally advised to use the CMA-ES \cite{hansen2003reducing,hansen2014principled,hansen2001completely,hansen2004evaluating,akimoto2020diagonal} as the oracle $M_{\epsilon_x,r_x}^{A_x}$ in AES instead of the (1+1)-ES.
We would like to emphasize that the difficulty addressed in this paper is associated with the interaction term between $x$ and $y$ in $f$, i.e., $H_{x,y}$ and $H_{y,x}$, and the CMA-ES does not help with this. Using the CMA-ES helps to solve poor conditioning of $H_{x,x}$ and $H_{y,y}$, which is not considered in this study.

As an important application of our approach, we are interested in a robust design of a feedback controller for automatic berthing.
When we design a controller, we first model the state equation of a control target, which is a ship in this case, and optimize the parameters of the controller on simulation.
However, there always exist modelling error and uncertainty in environment. For example, coefficients of a state equation model are typically estimated by water tank tests and weather conditions are uncertain when modelling. 
To design a controller reliable in a real environment, a task is formulated as a simulation-based min--max optimization problem \eqref{eq:spo}.
We apply our approach to robust automatic berthing to demonstrate the usefulness of the proposed approach. 
The results will be reported in another occasion.

We close the paper with open questions on this topic.
(I) A sufficient condition on $\eta$ derived in \Cref{thm:glc} for global linear convergence is not tight for $\delta > 0$ (see \Cref{sec:exp}).
We conjecture that we may be able to obtain a condition similar to that obtained for a simultaneous gradient update \cite{Mescheder2017nips}.
(II) The current analysis does not guarantee that our approach does not converge to nonsaddle critical points, whereas our experiments in \Cref{sec:exp} demonstrated that it does not.
(III) The proposed learning rate adaptation is yet to be theoretically analyzed for its convergence.
(IV) This study did not investigate methods for improving the convergence rate $\gamma$ when the interaction terms $H_{x,y}$ and $H_{y,x}$ are ill conditioned.
The proposed approach introduced a learning rate $\eta < 1$ to guarantee convergence, but the resulting convergence rate depends on $\eta$.
(V) There can exist a local (or global) optimum to the min--max optimization problem \eqref{eq:spo} that is not a min--max saddle point.
This corresponds to the case where the optimum is not a critical point of any of $y$.
Though it is out of the scope of the current paper, developing an algorithm that can locate such a local optimum of \eqref{eq:spo} is an important direction of the future work.

\begin{acks}
The author would like to thank the authors of \cite{Al-Dujaili2019lego} for publishing their code. 
This work is partially supported by JSPS KAKENHI Grant Number 19H04179. 
\end{acks}


\balance

\clearpage
\appendix
\section{Proofs}\label{apdx:proof}

The following proposition is derived immediately from the implicit function theorem, Theorem~5 of \cite{deoliveira2013}.
The proofs of \Cref{thm:glc,thm:llc} are based on this result.
\begin{proposition}[Implicit Function Theorem]\label{prop:implicit}
  In the following, let $(x^*, y^*)$ be a strict saddle point of $f$ and $f$ be locally strong convex--concave around $(x^*, y^*)$.
  There exist open sets $\mathcal{E}_{x,x} \subseteq \X$ including $x^*$ and $\mathcal{E}_{x,y} \subseteq \Y$ including $y^*$,
  such that there is a unique $\opty: \mathcal{E}_{x,x} \to \mathcal{E}_{x,y}$ such that $\nabla_y f(x, \opty(x)) = 0$. Moreover, $y^* = \opty(x^*)$ and $J_{\opty}(x) = -(H_{y,y}(x, \opty(x)))^{-1} H_{y,x}(x, \opty(x))$ for all $x \in \mathcal{E}_{x,x}$.
  Analogously, 
  there exist open sets $\mathcal{E}_{y,y} \subseteq \Y$ including $y^*$ and $\mathcal{E}_{y,x} \subseteq \X$ including $x^*$,
  such that there is a unique $\optx: \mathcal{E}_{y,y} \to \mathcal{E}_{y,x}$ such that $\nabla_x f(\optx(y), y) = 0$. Moreover, $x^* = \optx(y^*)$ and $J_{\optx}(y) = -(H_{x,x}(\optx(y), y))^{-1} H_{x,y}(\optx(y), y)$ for all $y \in \mathcal{E}_{y,y}$.
  If $f$ is globally strong convex--concave, one can take $\mathcal{E}_{x,x} = \mathcal{E}_{y,x} = \X$ and $\mathcal{E}_{y,y} = \mathcal{E}_{x,y} = \Y$ on the above statements. 
\end{proposition}

\Cref{prop:implicit} shows that under the condition stated in the proposition, the local minimal solutions to $\argmin_{x' \in \X} f(x', y)$ and $\argmin_{y' \in \Y} - f(x, y')$ are uniquely determined by $\optx(y)$ and $\opty(x)$. Moreover, by the mean value theorem we obtain
\begin{equation}
  \optx(y) - x^* = - (H_{x,x}(\optx(\tilde{y}), \tilde{y}))^{-1} H_{x,y}(\optx(\tilde{y}), \tilde{y}) (y - y^*) \enspace,
\label{eq:optx-opty}
\end{equation}
where $\tilde{y} = (1 - b) y + b y^*$ for some $b \in (0, 1)$. 
The approximate minimization oracle returns its approximate solutions $x'$ and $y'$. 
Then, we can write
\begin{multline}
  x_{t+1} - x^* = (1 - \eta)(x_t - x^*) \\
  + \eta (x_t' - \optx(y_t)) + \eta(\optx(y_t) - x^*) \enspace.
\label{eq:decompose}
\end{multline}
Analogously we can obtain the decomposition of $y_{t+1} - y^*$. 

\subsection{Proof of \Cref{thm:glc}}

  Let $x' = M^{A_x}_{\epsilon_x,r_x}(f(\cdot, y_t), x_t)$ and $y' = M^{A_y}_{\epsilon_y,r_y}(-f(x_t, \cdot), y_t)$.
  Let $\optx(y): \Y \to \X$ and $\opty(x): \X \to \Y$ be as shown in \Cref{prop:implicit}.
  We are going to investigate each term of \eqref{eq:decompose}.

  First, we study $x' - \optx(y_t)$ and $y' - \opty(x_t)$, which are the deviations of the outputs of the approximate minimization oracles from the exact solutions.
  For $(x_t, y_t) \in \X\times\Y$, we can write $\sqrt{\smash[b]{G_{x,x}^*}}(x' - \optx(y_t)) = E_x \sqrt{\smash[b]{G_{x,x}^*}} (x_t - \optx(y_t))$ and $\sqrt{\smash[b]{G_{y,y}^*}}(y' - \opty(x_t)) = E_y \sqrt{\smash[b]{G_{y,y}^*}}(y_t - \opty(x_t))$,
  where $E_x$ and $E_y$ are square matrices defined as
  \begin{equation*}
    \begin{split}
      E_{x} &= \sqrt{\smash[b]{G_{x,x}^*}}(x' - \optx(y_t))(x_t - \optx(y_t))^\T \sqrt{\smash[b]{G_{x,x}^*}} / \norm{x_t - \optx(y_t)}_{G_{x,x}^*}^2\\
      E_{y} &= \sqrt{\smash[b]{G_{y,y}^*}}(y' - \opty(x_t))(y_t - \opty(x_t))^\T \sqrt{\smash[b]{G_{y,y}^*}} / \norm{y_t - \opty(x_t)}_{G_{y,y}^*}^2\enspace.
    \end{split}
  \end{equation*}
  From \eqref{eq:oracle}, we have
  \begin{equation}
    \begin{split}
      &\textstyle \norm{x' - \optx(y_t)}_{G_{x,x}^*}^2 \leq \epsilon_x \norm{x_t - \optx(y_t)}_{G_{x,x}^*}^2 \enspace,
      \\
      &\textstyle \norm{y' - \opty(x_t)}_{G_{y,y}^*}^2 \leq \epsilon_y \norm{y_t - \opty(x_t)}_{G_{y,y}^*}^2 \enspace.
    \end{split}\label{eq:div-kappa}
\end{equation}
Then, from \eqref{eq:algo} and \eqref{eq:div-kappa}, we have
\begin{equation*}
  \begin{split}
    \sigma(E_x) &= \norm{x' - \optx(y_t)}_{G_{x,x}^*}^2 / \norm{x_t - \optx(y_t)}_{G_{x,x}^*}^2 \leq \epsilon_x \\
    \sigma(E_y) &= \norm{y' - \opty(x_t)}_{G_{y,y}^*}^2 / \norm{y_t - \opty(x_t)}_{G_{y,y}^*}^2 \leq \epsilon_y \enspace.
  \end{split}
\end{equation*}

  Next, we study $\optx(y_t) - x^*$ and $\opty(x_t) - y^*$, which are the differences between the saddle point $(x^*, y^*)$ and the solutions to the minimization problems that the approximate minimization oracles are solving. Let
  \begin{align*}
    \Delta_{x,y} =\textstyle (H_{x,x}(\optx(\tilde{y}), \tilde{y}))^{-1} H_{x,y}(\optx(\tilde{y}), \tilde{y})  - (H_{x,x}^*)^{-1} H_{x,y}^* ,
  \end{align*}
  where $\tilde{y}$ is as defined below \eqref{eq:optx-opty}.
  From \eqref{eq:optx-opty}, we can write
  \begin{align*}
    \optx(y_t)
    - x^*  = - ((H_{x,x}^*)^{-1} H_{x,y}^* + \Delta_{x,y}) \cdot (y_t - y^*) \enspace.
  \end{align*}
  Then we obtain
  \begin{multline*}
   \sqrt{\smash[b]{G_{x,x}^*}}(\optx(y_t) - x^*) = - \Big(\sqrt{\smash[b]{G_{x,x}^*}} (H_{x,x}^*)^{-1}H_{x,y}^* \sqrt{\smash[b]{G_{y,y}^*}}^{-1} \\
     + \sqrt{\smash[b]{G_{x,x}^*}} \Delta_{x,y} \sqrt{\smash[b]{G_{y,y}^*}}^{-1}\Big) \cdot \sqrt{\smash[b]{G_{y,y}^*}}(y_t - y^*) \enspace.
 \end{multline*}
 Analogously, in light of \Cref{prop:implicit}, we have
 \begin{align*}
    \opty(x_t)
    - y^* = - ((H_{y,y}^*)^{-1} H_{y,x}^* + \Delta_{y,x}) \cdot (x_t - x^*) \enspace,
 \end{align*}
 where
  \begin{align*}
    \Delta_{y,x} &=\textstyle (H_{y,y}(\opty(\tilde{x}), \tilde{x}))^{-1} H_{y,x}(\opty(\tilde{x}), \tilde{x}) - (H_{y,y}^*)^{-1} H_{y,x}^* , 
  \end{align*}
 where $\tilde{x} = (1 - a) x + a x^*$ for some $a \in (0, 1)$.  
 We obtain
  \begin{multline*}
    \sqrt{\smash[b]{G_{y,y}^*}}(\opty(x_t) - y^*) = - \Big(\sqrt{\smash[b]{G_{y,y}^*}} (H_{y,y}^*)^{-1}H_{y,x}^* \sqrt{\smash[b]{G_{x,x}^*}}^{-1} \\
     + \sqrt{\smash[b]{G_{y,y}^*}}\Delta_{y,x} \sqrt{\smash[b]{G_{x,x}^*}}^{-1} \Big) \cdot \sqrt{\smash[b]{G_{x,x}^*}}(x_t - x^*) \enspace.
  \end{multline*}
  
  Because of the decomposition \eqref{eq:decompose} and the same decomposition for $y$, i.e., 
  \begin{multline}
    y_{t+1} - y^* = (1 - \eta)(y_t - y^*) \\
    + \eta (y_t' - \opty(x_t)) + \eta(\opty(x_t) - y^*) \enspace,
  \end{multline}
  we have
  \begin{equation}
    \begin{bmatrix}
      \sqrt{\smash[b]{G_{x,x}^*}}(x_{t+1} - x^*) \\
      \sqrt{\smash[b]{G_{y,y}^*}}(y_{t+1} - y^*)
    \end{bmatrix}
 = (F + \eta \tilde{R}_1 + \eta\tilde{R}_2)
   \begin{bmatrix}
     \sqrt{\smash[b]{G_{x,x}}} (x_t - x^*)\\
     \sqrt{\smash[b]{G_{y,y}}} (y_t - y^*)   
   \end{bmatrix},
    \label{eq:matrixform}
  \end{equation}
  where
  \begin{align*}
    F &=
    \begin{bmatrix}
    (1-\eta) I  \quad - \eta \sqrt{\smash[b]{G_{x,x}^*}} (H_{x,x}^*)^{-1}H_{x,y}^* \sqrt{\smash[b]{G_{y,y}^*}}^{-1} \\
    - \eta \sqrt{\smash[b]{G_{y,y}^*}} (H_{y,y}^*)^{-1}H_{y,x}^* \sqrt{\smash[b]{G_{x,x}^*}}^{-1} \quad (1-\eta) I
      \end{bmatrix},
    \\
     \tilde{R}_1 &=
    \begin{bmatrix}
    E_x & 0 \\ 0 & E_y
  \end{bmatrix},\                    
    \tilde{R}_2 =
   \begin{bmatrix}
    0 \quad - \sqrt{\smash[b]{G_{x,x}^*}}\Delta_{x,y} \sqrt{\smash[b]{G_{y,y}^*}}^{-1} \\
    -  \sqrt{\smash[b]{G_{y,y}^*}}\Delta_{y,x} \sqrt{\smash[b]{G_{x,x}^*}}^{-1} \quad 0
   \end{bmatrix}.
\end{align*}
  
Now, we show that there exists $\gamma > 0$ such that for any $(x_t, y_t) \in \X\times\Y$, 
  \begin{multline}
   \norm{x_{t+1} - x^*}_{G_{x,x}^*}^2 +  \norm{y_{t+1} - y^*}_{G_{y,y}^*}^2\\
    \leq \gamma^2 \big(\norm{x_{t} - x^*}_{G_{x,x}^*}^2 +  \norm{y_{t} - y^*}_{G_{y,y}^*}^2\big)
 \enspace.
    \label{eq:llc-sub}
  \end{multline}
  Notably, the left-hand side of \eqref{eq:llc-sub} is the squared Euclidean norm of the left-hand side of \eqref{eq:matrixform}.
  Therefore, to prove \eqref{eq:llc-sub}, it suffices to show that the greatest singular value of $F + \eta \tilde{R}_1 + \eta \tilde{R}_2$ has an upper bound $\gamma$. 
The greatest singular value of $F + \eta \tilde{R}_1 + \eta \tilde{R}_2$ has an upper bound of the sum of the greatest singular value of each term, i.e., $\sigma(F + \eta \tilde{R}_1 + \eta \tilde{R}_2) \leq \sigma(F) + \eta \sigma(\tilde{R}_1) + \eta \sigma(\tilde{R}_2)$. 
The greatest singular value of $F$ is equal to the square root of the greatest eigenvalue of $F F^\T$ in this case.
Noting that $\sqrt{\smash[b]{G_{x,x}}} (H_{x,x}^*)^{-1}H_{x,y}^* \sqrt{\smash[b]{G_{y,y}}}^{-1} = \sqrt{\smash[b]{G_{x,x}}}^{-1} H_{x,y}^* (-H_{y,y}^*)^{-1} \sqrt{\smash[b]{G_{y,y}}}$, 
we have
\begin{multline*}
  FF^\T
  = (1 - \eta)^2 I +
  \eta^2 \times \\
    \begin{bmatrix}
      \sqrt{\smash[b]{G_{x,x}}} (H_{x,x}^*)^{-1} H_{x,y}^* G_{y,y}^{-1} H_{y,x}^* (H_{x,x}^*)^{-1} \sqrt{\smash[b]{G_{x,x}}} \qquad 0  \\
      0 \qquad \sqrt{\smash[b]{G_{y,y}}} (H_{y,y}^*)^{-1}H_{y,x}^* G_{x,x}^{-1} H_{x,y}^* (H_{y,y}^*)^{-1} \sqrt{\smash[b]{G_{y,y}}}
    \end{bmatrix} \enspace,
\end{multline*}
whose greatest eigenvalue is $(1-\eta)^2 + \eta^2 \bar\sigma^2$. Hence, $\sigma(F) = ((1-\eta)^2 + \eta^2 \bar\sigma^2)^{1/2}$. 
The greatest singular value of $\tilde{R}_1$ has an upper bound of $\bar\epsilon = \max(\epsilon_x, \epsilon_y)$, as derived below \eqref{eq:div-kappa}. The greatest singular value of $\tilde{R}_2$ has an upper bound of
$$\max\big[ \sigma\big(\sqrt{\smash[b]{G_{x,x}^*}}\Delta_{x,y}\sqrt{\smash[b]{G_{y,y}^*}}^{-1} \big), \sigma\big(\sqrt{\smash[b]{G_{y,y}^*}}\Delta_{y,x}\sqrt{\smash[b]{G_{x,x}^*}}^{-1} \big) \big]\leq \delta\enspace.$$ 
  Therefore, $\gamma \leq ((1-\eta)^2 + \eta^2 \bar\sigma^2)^{1/2} + \eta (\bar\epsilon + \delta) < 1$.

Because \eqref{eq:llc-sub} holds for any $(x_t, y_t) \in \X\times\Y$, a repeated application of \eqref{eq:llc-sub} shows \eqref{eq:glc-conclusion}. This completes the proof.

\subsection{Proof of \Cref{thm:llc}}

  First, we show that there exists $\bar{r}_x > 0$ and $\bar{r}_y > 0$ such that 
  for any $r_x \in (0, \bar{r}_x)$ and $r_y \in (0, \bar{r}_y)$, there exists $U_x$ and $U_y$ such that 
  $\opty(x) = \argmin_{y' \in U_{A_y}(y,r_y)} -f(x, y')$ and $\optx(y) = \argmin_{x' \in U_{A_x}(x,r_x)} f(x', y)$ hold for all $(x, y) \in U_x \times U_y$,
  where $\opty: U_x \to U_{A_y}(y,r_y)$ and $\optx: U_y \to U_{A_x}(x,r_x)$ are as defined in \Cref{prop:implicit}.
  In light of \Cref{prop:implicit}, there exists a unique $\optx: \mathcal{E}_{x,x} \to \mathcal{E}_{x,y}$ and a unique $\opty: \mathcal{E}_{y,y} \to \mathcal{E}_{y,x}$
  for some $\mathcal{E}_{x,x}$, $\mathcal{E}_{x,y}$, $\mathcal{E}_{y,x}$, and $\mathcal{E}_{y,y}$.
  Because $A_x = G_{x,x}^{*}$, we have $U_{A_x}(x, r_x) \subseteq \{x' \in \X : \norm{x' - x}_{G_{x,x}^{*}} \leq r_x\}$.
  Hence, there exists a $\bar{r}_x > 0$ such that for all $r_x \in (0, \bar{r}_x)$ we have $U_{A_x}(x^*, r_x) \subseteq \mathcal{E}_{y,x}$.
  Analogously, we have that there exists a $\bar{r}_y > 0$ such that for all $r_y \in (0, \bar{r}_y)$ we have $U_{A_y}(y^*, r_y) \subseteq \mathcal{E}_{x,y}$.
  Because $\optx$ and $\opty$ are continuous and $\optx(y^*) = x^*$ and $\opty(x^*) = y^*$,
  for any $r_x \in (0, \bar{r}_x)$ and $r_y \in (0, \bar{r}_y)$, there exist $U_x \subseteq \mathcal{E}_{x,x}$ including $x^*$ and $U_y \subseteq \mathcal{E}_{y,y}$ including $y^*$ such that $\optx(y) \in U_{A_x}(x, r_x)$ and $\opty(x) \in U_{A_y}(y, r_y)$ for all $x \in U_x$ and $y \in U_y$. 
  
  For $(x, y) \in U_x\times U_y$, we obtain \eqref{eq:llc-sub} by following the line of the proof of \Cref{thm:llc}, where
  $\gamma \leq ((1-\eta)^2 + \eta^2 \bar\sigma^2)^{1/2} + \eta (\bar\epsilon + \delta)$ with $\delta$ being the maximum of $\sigma\big(\sqrt{\smash[b]{G_{x,x}^*}}\Delta_{x,y}\sqrt{\smash[b]{G_{y,y}^*}}^{-1}\big)$ and $\sigma\big(\sqrt{\smash[b]{G_{y,y}^*}}\Delta_{y,x}\sqrt{\smash[b]{G_{x,x}^*}}^{-1}\big)$.
  Because of the continuities of $H_{x,x}$, $H_{x,y}$, $H_{y,x}$, and $H_{y,y}$,
  for any $\bar\delta > 0$ there exists a neighborhood $\tilde{U} \subseteq U_x \times U_y$ of $(x^*, y^*)$ such that
  $\bar\delta < \bar\delta$ for all $(x, y) \in \tilde{U}$.
  That is, for any $\gamma \in (\bar\gamma, 1)$, where $\bar\gamma = ((1-\eta)^2 + \eta^2 \bar\sigma^2)^{1/2}
  + \eta \bar\epsilon$, we can find a neighborhood $\tilde{U}$ of $(x^*, y^*)$ such that \eqref{eq:llc-sub} holds for all $(x_t, y_t) \in \tilde{U}$. 

  Finally, we show that there exists $U \subseteq \tilde{U}$ such that
  the sequence $(x_t, y_t)$ starting from $(x_0, y_0) \in U$ never leaves $U$.
  Define $B(r) := \{(x, y) \in \X\times\Y : \tilde{G}(x, y) < r^2/2\}$.
  Then, one can find an $\bar{r} > 0$ such that $B(\bar{r}) \subseteq \tilde{U}$.
  Let $U = B(\bar{r}) \subseteq \tilde{U}$.
  Because $\gamma < 1$, \eqref{eq:llc-sub} implies that if $(x_0, y_0) \in U$, we have $(x_t, y_t) \in U$.
  Hence, $(x_t, y_t)$ never leaves $U$.

  Repeated applications of \eqref{eq:llc-sub} lead to \eqref{eq:llc-conclusion}. This completes the proof.

\end{document}